\documentclass[12pt]{article} 
\usepackage{amssymb} 
 
\usepackage{psfrag, graphicx} 
\usepackage{amsmath} 
\usepackage{amsthm}
 
\theoremstyle{plain} 
\newtheorem{theorem}{Theorem}[section] 
\newtheorem{thm}[theorem]{Theorem} 
\newtheorem{lemma}[theorem]{Lemma} 
\newtheorem{cor}[theorem]{Corollary} 
\newtheorem{corollary}[theorem]{Corollary} 
\newtheorem{prop}[theorem]{Proposition}

\theoremstyle{definition}
\newtheorem{defn}[theorem]{Definition} 
\newtheorem{conj}[theorem]{Conjecture} 
 
\newtheorem{example}[theorem]{Example} 
\theoremstyle{remark}
\newtheorem{rmk}[theorem]{Remark} 
\newtheorem{remark}[theorem]{Remark} 
\newtheorem{obs}[theorem]{Observation}

\newcommand{\be}{\begin{equation}} 
\newcommand{\ee}{\end{equation}} 
\newcommand{\ben}{\begin{enumerate}} 
\newcommand{\een}{\end{enumerate}} 
\newcommand{\bi}{\begin{itemize}} 
\newcommand{\ei}{\end{itemize}} 
\newcommand{\ra}{{\rightarrow}} 

\providecommand{\abs}[1]{\lvert#1\rvert}
\providecommand{\norm}[1]{\lVert#1\rVert} 

\newcommand{\wtilde}{\widetilde} 

\newcommand{\bC}{{\mathbb C}} 
 
\newcommand{\bR}{{\mathbb R}} 
\newcommand{\bZ}{{\mathbb Z}} 
\newcommand{\cA}{{\cal A}}

\newcommand{\cS}{{\cal S}} 
\newcommand{\bdy}{\partial} 
\newcommand{\del}{\partial}

\newcommand{\atf}{almost toric fibration\ }
\newcommand{\triple}{$(B,\cA,\cS)\ $}

\begin{document} 
 
\title{Almost toric symplectic four-manifolds}

\author{ 
Naichung Conan Leung\footnote{
School of Mathematics, University of Minnesota, MN 55455; leung@math.umn.edu
 \newline
        Supported in part by NSF/DMS-0103355}
 \ and \ Margaret Symington\footnote{
School of Mathematics, Georgia Inst. of Tech., Atlanta, GA
30332; msyming@math.gatech.edu
\newline
Supported in part by NSF/DMS-0204368}}

\maketitle
 
\abstract{
Almost toric manifolds form a class of singular Lagrangian fibered
symplectic manifolds that is a natural generalization of toric manifolds.
Notable examples include the K3 surface, the phase space of the 
spherical pendulum and rational balls useful for symplectic surgeries. 
The main result of the paper is a complete classification up to diffeomorphism
of closed almost toric four-manifolds.

A key step in the proof is a geometric classification of the singular affine
structures that can occur on the base of a closed almost toric four-manifold.}

\section{Introduction}

Almost toric manifolds, introduced by the second author 
in~\cite{Symington.fourfromtwo}, are symplectic manifolds equipped with
a fibration structure that generalizes 
toric manifolds while retaining some of their geometric features and rigidity.
Accordingly, almost toric manifolds lie at the interface of
symplectic topology, toric geometry, integrable systems and, in
dimension four, mirror symmetry.
They enjoy the property (similar to toric manifolds)
that much symplectic and topological information is encoded in the base of
the fibration;
they can be used to efficiently describe certain symplectic surgeries such
as symplectic sums and rational blowdowns~\cite{Symington.grbd};  
they accommodate singularities that are typical in an integrable system
(focus-focus singularities);
furthermore, generic special Lagrangian fibrations of K3 surfaces -- 
of interest in mirror symmetry (cf.~\cite{StromingerYauZaslow.mirror},~\cite{Morrison.geomofmirrorsymm},~\cite{GrossWilson.largecxlimit} 
and~\cite{GrossSiebert.affinelog}) -- 
are almost toric fibrations.

In dimension four, an informal definition of an almost toric manifold is a 
symplectic four-manifold $(M,\omega)$ equipped with a projection 
$\pi:(M,\omega)\ra B$ 
to a surface $B$ such that locally $\pi$ has the structure of either the
moment map for a torus action or the Lagrangian analog of a
Lefschetz fibration.  
The main result of this paper (Theorem~\ref{diffclass.thm})
is the classification up to diffeomorphism of closed four-manifolds 
with such a structure.  The theorem, which also includes information
about the structure of possible fibrations, is summarized in 
Table~\ref{diffclass.table}.

\begin{table} 
\begin{tabular}{|c|c|c|c|} 

\hline
& & & \\

Base & \# of nodes & \# of vertices & Total space \\
& & &  \\ 

 \hline  & & & \\
 \hline & & &  \\

$D^2$ & $n \ge 0$ & $k \ge {\rm max}(0,3-n) $ & $\bC P^2 \# (n+k-3)\overline{\bC P}^2$ or \\ 
 & & & $ S^2\times S^2$ (if $n+k=4$)  \\ 
& & & \\ \hline & & &  \\

$S^1\times I$ & $n\ge 0$ & $0$ & $(S^2\times T^2)\,\#\, n\overline{\bC P}^2$ 
or \\ 
& & & $(S^2\tilde \times T^2)\,\#\, n\overline{\bC P}^2$ \\ & & &  \\ 
 \hline & & & \\

$S^1\tilde \times I$ & $n\ge 0$ & $0$ & $(S^2\times T^2)\,\#\, n\overline{\bC P}^2$ 
or \\ 
& & & $(S^2\tilde \times T^2)\,\#\, n\overline{\bC P}^2$ \\ & & &  \\ 
\hline & & & \\

$S^2$ & $24$ & $0$ & K3 surface \\ & & &  \\ \hline & & & \\

$\bR P^2$ & $12$ & $0$ & Enriques surface \\ & & & \\ \hline & & & \\

$T^2$ & $0$ & $0$ & $T^2$ bundle with \\
& & & monodromy  
$ 
\left(I,\left(\begin{smallmatrix}1&\lambda\\0&1\end{smallmatrix}\right)\right)
$ \\
& & &  \\ \hline & & & \\

$\ S^1\tilde \times S^1\ $ & $0$ & $0$ & $\ T^2$ bundle with \ \\
& & & monodromy  
$\left(\left(\begin{smallmatrix}1&0\\0&-1\end{smallmatrix}\right),
\left(\begin{smallmatrix}1&\lambda\\0&1\end{smallmatrix}\right)\right)$\\
& & &  \\ \hline 
\end{tabular}
\caption{Closed almost toric four-manifolds}
\label{diffclass.table} 
\end{table}

The list of closed four-manifolds admitting an almost toric fibration is 
fairly short due to the Lagrangian constraint.  
For instance, the only manifolds that admit a fibration that is locally \lq\lq
Lagrangian Lefschetz\rq\rq\ 
are the K3 surface and a $\bZ_2$ quotient of the K3 that is diffeomorphic to
the Enriques surface.  

It is straightforward, except for the Enriques surface,  
to deduce from the existing literature 
(\cite{Zung.II}, \cite{Geiges.torusbundles},
\cite{SakamotoFukuhara.torusbundles}) that each of the manifolds
listed in Table~\ref{diffclass.table}
does indeed admit an almost toric fibration.  
(While one would certainly expect the $\bZ_2$ quotient of the 
K3 to be the Enriques surface,  this requires proof since the base is 
$\bR P^2$ rather than $\bC P^1$ as in the holomorphic case; 
see Lemma~\ref{Enriques.lem}).
In order to show that this is the complete list we appeal
to the work of Zung~\cite{Zung.II} which shows the degree to which the
geometry of the base of an almost toric fibration controls the topology
of the total space.  The primary task then is to study the possible
structures on the base of an almost toric fibration.  

{\sc Acknowledgments:}  The second author would like to thank John Etnyre and
Eugene Lerman for pointing out useful references, Bill Goldman for
providing notes on an argument of Benzecri~\cite{Benzecri.affine} and
Ilia Zharkov for comments on a preliminary draft.

\section{Background and results} \label{background.sec}
 
A symplectic manifold of dimension $2n$ is toric if it is 
equipped with an effective Hamiltonian $T^n$ action.
Toric manifolds are well studied for their beautiful geometric 
properties (e.g.~\cite{Audin.torus}, \cite{Guillemin.toric})
and their relevance to mirror symmetry (cf.~\cite{Batyrev.mirrortoric}).
Delzant's Theorem (\cite{Delzant.moment}) asserts a fundamental property of 
closed toric manifolds: the manifold, 
symplectic structure and torus action are completely determined by a 
polytope in $\bR^n$, the image of the {\it moment map}.
The preimages of regular values of the moment map $\pi:(M,\omega)\ra\bR^n$
are Lagrangian tori (half-dimensional tori on which the symplectic form 
vanishes).
Furthermore, any critical point of $\pi$ is an {\it elliptic} singular
point: it has a Darboux neighborhood (with symplectic form 
$dx\wedge dy:=\sum_i dx_i\wedge dy_i$) 
in which the map $\pi:=(\pi_1,\ldots \pi_k,\pi_{k+1},\ldots \pi_{n})$
has the form $\pi_j(x,y)=x_j$ for $j\le k$ and 
$\pi_j (x,y)=(x_{j}^{2}+y_{j}^{2})$ for $j>k$.
This provides the base of the fibration with a stratification (which
we denote by $\cS$) according
to the dimension of its preimage.  In the case of a closed manifold this
stratification agrees with the natural stratification of the polytope
image according to the dimension of the facets.
Furthermore, the preimage of singular values are tori that
are Lagrangian submanifolds of the preimage of the stratum containing the
singular value.

Allowing, as in the holomorphic case, for a fibration to have singular 
fibers, a toric manifold provides an example of a Lagrangian fibration:
\begin{defn}
A {\it Lagrangian fibration} of a symplectic manifold $(M,\omega)$
(possibly with boundary)
is a map $\pi :(M,\omega )\rightarrow B$ to a space of half the dimension 
such that on the preimage of an
open dense set $B_{0}\subset B$ the projection $\pi$ is a locally trivial
fibration with $\omega |_{\pi^{-1}(b)}=0$ for all $b\in B_0$. 
We assume that the fibers of a Lagrangian fibration are
compact, connected and without boundary.
\end{defn}

Noting that the moment map for a toric manifold provides an immersion of the
base to 
$\bR^n$, there are two natural ways to generalize toric manifolds within 
the class of Lagrangian fibrations:  
one can allow the base to be a space that does not immerse in $\bR^n$ and
one can allow more general singular fibers (than lower dimensional tori).
Boucetta and Molino~\cite{BoucettaMolino.fibrations} have considered the 
first type of generalization, establishing complete invariants.
Zung~\cite{Zung.II} 
generalized still further, defining a rather general class of 
Lagrangian fibrations -- in particular those with non-degenerate
topologically stable singularities -- and specifying 
the data required to classify such
fibrations up to fiber-preserving symplectomorphism.
Almost toric fibrations form a subset of these in which we exclude
hyperbolic singularities:

\begin{defn}  \label{almosttoric.defn}
An {\it almost toric fibration} of a symplectic $2n$-manifold is
a Lagrangian fibration $\pi:(M,\omega)\ra B$ such that 
any critical point of $\pi$ has a Darboux neighborhood (with symplectic 
form $dx\wedge dy$) 
in which the projection $\pi:=(\pi_1,\ldots \pi_k,\pi_{k+1},\ldots \pi_{n})$
where $\pi_j(x,y)=x_j$ for $j\le k$ and the other components
have one of the following two forms:

\begin{align} 
\pi_j (x,y)&=(x_{j}^{2}+y_{j}^{2}) \quad \textit{elliptic} {\rm \ or\ }
\textit{toric}  \\
(\pi_i,\pi_j) (x,y)&= (x_{i}y_{i}+x_{j}y_{j},
         x_{i}y_{j}-x_{j}y_{i}) \quad \textit{nodal, 
{\rm or} focus-focus}.
\end{align}
An {\it almost toric manifold} is a symplectic manifold equipped with
an almost toric fibration.
A {\it toric fibration} is a Lagrangian fibration induced by an effective
Hamiltonian torus action.
\end{defn}

An alternative way to present almost toric manifolds is as manifolds
with the local structure of an integrable system having compact 
fibers and only elliptic and nodal singularities (or products thereof)
where the nodal singularities cause only positive intersections.  
Recall that an {\it integrable system} is a symplectic 
$2n$-manifold $(M,\omega)$
equipped with a collection of $n$ functionally independent
Poisson commuting functions $F_i:M\rightarrow \bR^n$.
For an almost toric manifold one can take the $F_i$ 
to be the components of $\Pi\circ\pi$ where $\Pi:U\ra \bR^n$ is a coordinate 
chart on the base in $U\subset B$ 

\begin{rmk} \label{positivity.rmk}
The self-intersection that appears in a fiber with a nodal (focus-focus) 
singularity is always positive~\cite{Zung.focusfocus}.  
While Lagrangian planes are not by themselves oriented, any orientation
of the base orients these planes (via the Hamiltonian vector fields induced
by a basis in the cotangent bundle of the base) thereby giving a well-defined
sign to the intersection.  
\end{rmk}

\begin{rmk} 
The reason for excluding 
hyperbolic singularities from the definition of almost toric manifolds 
is that they greatly complicate the process of
recovering the total space from the base.
This is evident already in dimension two where the preimage of a singular
value can include several hyperbolic singular points and the separatrices
that connect them.
\end{rmk}

Assigning points in the base of an almost toric fibration
to strata according to the dimension of
their preimage yields, as in the toric case, a stratification $\cS$ 
of the base~\cite{Zung.II}.
In dimension four the images of nodal singular points ({\it nodes}) 
are isolated points that belong to the top dimensional stratum.

The base of an almost toric fibration also carries a fairly rigid geometric 
structure defined on the regular values of the fibration map.

\begin{defn} An {\it integral affine structure} $\cA$ on a manifold $B$
(possibly with boundary) is a lattice in its tangent bundle.  
A manifold admitting such a structure is an {\it integral affine manifold}.
\end{defn}
Alternatively, one could define an integral affine $n$-manifold to be a
manifold whose structure group is ${\rm Aff}(n,\bZ)=GL(n,\bZ)\ltimes \bR^n$.

The integral affine structure on the base of a regular Lagrangian fibration
arises from a natural action of the cotangent bundle of the base on the
total space: any $\alpha\in T^*B$ defines a vertical vector field $X_\alpha$ 
symplectically dual to $\alpha$, so $\alpha\cdot x=\varphi_\alpha(x)$ 
where $\varphi_\alpha(x)$ is the time-one flow of $X_\alpha$ is the natural
action.
The elements of the cotangent bundle of $B$ that act trivially
form a lattice $\Lambda^*$.  The dual lattice in
the tangent bundle defines the integral affine structure.
(See~\cite{Duistermaat.actionangle} or Section 2 
of~\cite{Symington.fourfromtwo} for more details.)
For any $\alpha\in\Lambda^*$ we denote by 
$[\alpha]$ the homology class of integral curves of the time-one
flow of $X_\alpha$.

A singular point is called an {\it elliptic singularity of corank} 
$k$~$(\ge 1)$ if
in the normal form of Definition~\ref{almosttoric.defn} the projection
map $\pi$ has $k$ elliptic factors and no nodal factors.
The locus of elliptic singular points of a given corank form a submanifold 
whose image is an integral affine submanifold of dimension $n-k$ in
the boundary of the base.
Furthermore, the normal form for elliptic singularities assures that the
base is a manifold with corners.

Any almost toric fibration is a locally toric fibration over $B-\Sigma$ where
$\Sigma$ is the codimension two set of points containing a nodal
singularity in their preimage.  (For an almost toric four-manifold
$\Sigma$ is a finite set of points on the interior of the base.)
By the integral affine structure $\cA$ on the base
of an almost toric fibration we mean the affine structure defined 
on $B-\Sigma$.  Likewise, the smooth structure on the
base is understood to be the smooth structure induced by $\cA$ on the
complement of the nodes.

If an affine structure $\cA$ and stratification $\cS$ are induced from
a toric or almost toric fibration, we call the triple $(B,\cA,\cS)$ a
{\it toric base} or an {\it almost toric base} respectively.
In both cases $(B,\cA,\cS)$ is a strong
invariant that dimension four often determines the
the total space (Corollary~\ref{basedetermination.cor}). 

A fiber in an almost toric fibration of dimension four can have more than one
singular point, but the vanishing cycles for the different singular points
must represent the same homology class in a regular fiber.
This forces a singular fiber with $k\ge 2$ nodal singularities to be a
reducible fiber with $k$ irreducible components, each diffeomorphic to
a sphere.  Furthermore, the number of singular points can be detected from
the integral affine structure near the image of the fiber 
(Section~\ref{nodalnbhd.sec}).
Throughout the paper the term {\it nodal fiber} refers to a fiber with
just one nodal singular point in the fiber (the generic
case), i.e. the Lagrangian analog of a Lefschetz fiber.
Accordingly, unless we specify that a node has multiplicity,
we assume that the preimage of a node contains exactly one singular point.

Nodal fibers arise naturally in Lagrangian fibrations both in the 
algebraic and integrable systems settings.  In the complex algebraic setting
one can start with a holomorphically
fibered K3 surface with $24$ singular (Lefschetz) fibers -- the generic case--
and perform a hyperk\"ahler rotation to make that same fibration Lagrangian.
Meanwhile, physical integrable systems with two degrees of freedom often
contain nodal fibers.  One example is the spherical pendulum 
(Example~\ref{pendulum.ex} in Section~\ref{noncpct.sec}).

We are now ready to state our main theorem.

\begin{thm} \label{diffclass.thm}
If $\pi:(M,\omega)\ra (B,\cA,\cS)$ is an almost toric fibration of a closed
four-manifold then
the total space $M$ must be diffeomorphic to 

(i) $S^2\times S^2$, 

(ii) $S^2\times T^2$, 

(iii) $N\,\#\,n\overline{\bC P}^2$ with
$N=\bC P^2$ or $S^2\tilde\times T^2$ and $n\ge 0$, 

(iv) the K3 surface, 

(v) the Enriques surface,

(vi) a torus bundle over the torus with monodromy
\begin{equation}
\left\{I,\left(  
\begin{array}{cc} 
1 & \lambda \\  
0 & 1 
\end{array} 
\right)\right\}, \ \  \lambda\in\bZ,
\end{equation} 
or

(vii) a torus bundle over the Klein bottle with monodromy 
\begin{equation}
\left\{\left( 
\begin{array}{cc} 
1 & 0 \\  
0 & -1 
\end{array} 
\right), \left(\begin{array}{cc} 
1 & \lambda \\  
0 & 1 
\end{array} 
\right)\right\}, \ \ \lambda\in\bZ.
\end{equation}
Furthermore, Table~\ref{diffclass.table} classifies such fibrations according
to the homeomorphism type of the base $B$, the number of nodes in the
affine structure $\cA$ and the number of vertices on the boundary of 
$(B,\cA,\cS)$ (i.e. the cardinality of the zero-stratum of $\cS$).
\end{thm}
Note that the sum of the number of nodes and vertices in the base equals
the Euler characteristic of the total space.
See Example~\ref{toriovertori.ex} for more details on the diffeomorphism types
of torus bundles over tori.

\begin{rmk}
We reiterate that the point of Theorem~\ref{diffclass.thm} is the brevity of 
the list of closed four-manifolds admitting an almost toric fibration.
It is easy to show, and not surprising, that 
the list of Lagrangian fibrations with elliptic singularities
is restricted to toric manifolds, sphere bundles over the torus 
(fibering over the cylinder or M\"oebius band) and torus bundles
over the torus or Klein bottle.
Therefore the question is, how much flexibility is gained by the allowance
of nodal fibers?  Since one can turn a holomorphic fibration of the K3
surface into a Lagrangian fibration via a hyperk\"ahler rotation it is
immediate that the K3 surface and its $\bZ_2$ quotient fibering over 
$\bR P^2$ admit almost toric fibrations.  
The possibility of blowing up at points in the preimage
of the one-stratum, first observed by Zung~\cite{Zung.II}, implies the
existence of 
almost toric fibrations of blowups of the sphere bundles over tori.
With this perspective, Theorem~\ref{diffclass.thm} 
states that this is the extent of flexibility introduced by nodal fibers.
\end{rmk}

Zung's work on (singular) Lagrangian fibrations~\cite{Zung.II}, together 
with uniqueness of the local structure of a neighborhood of a node 
(Proposition~\ref{nodalgerm.prop}), implies that the
base $(B,\cA,\cS)$ determines the topology of the total space of an
almost toric fibration in dimension four, except in the case of fibrations 
over the torus or Klein bottle.
Since the topology of the total spaces of Lagrangian fibrations over the
torus was already known (\cite{Geiges.torusbundles}, \cite{Smith.torusfibr}),
the main steps in the proof of Theorem~\ref{diffclass.thm} are
first establishing that the base of the fibration must be homeomorphic
to one of the manifolds listed in Table~\ref{diffclass.table} and
second proving any manifold that admits an almost toric fibration over
the disk also admits a toric fibration over the disk.

The first main step is accomplished by making precise the heuristic
that nodal singularities contribute non-negative curvature to the base. 
(Note that the non-trivial monodromy around a node precludes the existence
of a metric compatible with the affine structure.)  
The second step relies on an 
iterative process that transforms an almost toric base into a toric base
while preserving the topology (but not the symplectic structure) of the
total space.

\section{Total space from the base}

\subsection{Toric manifolds}  \label{toric.sec}

Delzant's theorem asserts that for a closed toric manifold of
dimension $2n$ the image of the moment map, a polytope in $\bR^n$,
determines the total space, its symplectic structure and the torus action.
If one drops the assumption that the total space is a closed manifold
(allowing noncompactness and nonempty boundary)
then many different manifolds can have the same moment map image.

The ambiguity arises in trying to determine the topology of the preimage of 
a value of the moment map: while it must be comprised of tori, 
neither the number of components nor their dimensions can be decided from the 
moment map image.  
For instance, if $\Delta\subset\bR^n$ is the
polytope that is the moment map image of a closed manifold, then 
$\Delta\times T^n$, with coordinates $(p,q)$ and symplectic structure 
$dp\wedge dq$, also admits a torus
action with $\Delta$ as its moment map image.  (In the integrable
systems language these are action-angle coordinates.)
Also, given any affine $n$-manifold $(V,\cA)$ for which there is
a surjective immersion onto $(\Delta,\cA_0)$, $V\times T^n$ also has a toric
fibration with moment map image $\Delta$.
Ambiguities concerning the 
topology of fibers can be addressed by using the toric base $(B,\cA,\cS)$
rather than the moment map image.

Thus, in the noncompact case we have the following analog of Delzant's
theorem in which we can recover the total space and Lagrangian fibration
but not a torus action.

\begin{prop} (cf. \cite{Symington.fourfromtwo}) \label{toricbase.prop}
If $(B,\cA,\cS)$ is a toric base of dimension $n$ then there is unique
symplectic manifold $(M,\omega)$ that admits a unique
Lagrangian fibration $\pi:(M,\omega)\ra (B,\cA,\cS)$.  
\end{prop}

To understand the construction, first note that 
the integral affine structure $\cA$ is intimately associated with 
the moment map.
Indeed, if $(B,\cA,\cS)$ is a  toric base
and $\Delta\subset\bR^n$ is the image of the corresponding
moment map, then there is an immersion $\Phi: (B,\cA)\ra (\Delta,\cA_0)$. 
Furthermore, the image of any immersion of $(B,\cA)$ into $(\bR^n,\cA_0)$ 
differs from $\Phi(B)$ only by an element of ${\rm Aff}(2,\bZ)$.

Now proceed as follows:  choose an affine immersion 
$\Phi:(B,\cA)\ra (\bR^n,\cA_0)$ which provides local coordinates $p$ on
any neighborhood of $B$ that embeds in $\bR^n$ via $\Phi$.
Consider a toric fibration
$\pi': (B\times T^n, dp\wedge dq)\ra (B,\cA)$.
The points of $(B,\cA,\cS)$ that belong to each connected component of the
$(n-1)$-stratum comprise a portion of $\bdy B$ whose preimage is fibered
by circles that are in the kernel of the $\omega$ restricted to 
$\bdy B\times R^n$.  Collapsing these circles on the preimage of the
closure of each component yields $(M,\omega)$ 
and a toric fibration $\pi: (M,\omega)\ra (B,\cA,\cS)$. 
(It is the assumption that $(B,\cA,\cS)$ is a toric base that ensures the
resulting space is a manifold.  See Section $3$ 
of~\cite{Symington.fourfromtwo} for an intrinsic characterization of a toric
base.)

Note that there is a symplectic
projection $\rho:(B\times T^n, dp\wedge dq)\ra(M,\omega)$
that is a diffeomorphism over the points $x$ of $M$ such that $\pi(x)$ 
belongs to the top-dimensional stratum of $(B,\cA,\cS)$.
Indeed, this presentation gives local action-angle 
coordinates $(p,q)$ on a dense subset of $(M,\omega)$, with
$(\frac{\del}{\del p_1}, \ldots \frac{\del}{\del p_n})$ being a basis for
the lattice in the tangent space at any point of $(B,\cA)$.
Furthermore, for some choice of $\Phi$ the original torus action is given on
the preimage of regular values of $\pi$ by $t\cdot (p,q)=(p,q+t)$, and this 
free action extends uniquely to all of $M$.

The procedure of constructing $(M,\omega)$ from $(B\times T^n,dp\wedge dq)$
has been called {\it boundary reduction} by the second 
author~\cite{Symington.fourfromtwo}.  

\subsection{Affine structure and monodromy} \label{monodromy.sec}

As we shall see in Section~\ref{atbase.sec},  
while Proposition~\ref{toricbase.prop} does not generalize
completely to the almost toric case, it comes close.  In
many cases the base of an almost toric fibration does determine
the total space.

An essential way in which the base $(B,\cA,\cS)$ 
influences the topology of the total space is by capturing the
monodromy.  Specifically, the topological monodromy of
the torus fibration over the regular values $B_0\subset B$
is determined by the affine monodromy in the base, i.e. the lattice in 
$T B_0$ (or, dually, in $T^* B_0$). 

The {\it affine monodromy} of an integral affine manifold $B$ is defined 
analogously to the monodromy of a torus fibration 
(cf.~\cite{GompfStipsicz.4mflds}).
Specifically, if $\Lambda$ is the lattice in $TB$, choose a point $b\in B_0$,
identify $(T_bB,\Lambda_b)$ with $(\bR^n,\bZ^n)$ and
for each element $\alpha\in \pi_1(B_0,b)$ choose a representative 
$\gamma_\alpha:I\ra B_0$.
The monodromy representation relative to
these choices is $\Psi_B:\pi_1(B,b)\rightarrow {\rm Aff}(n,\bZ)$ where
$\Psi_B(\alpha)$ is the automorphism of $(\bR^n,\bZ^n)$ such that 
$\gamma_\alpha^* (TB,\Lambda)$ is isomorphic 
to $I\times (\bR^n,\bZ^n)/((0,p)\sim (1,\Psi_B(\alpha)(p))$, $p\in\bR^n$.
The monodromy is the equivalence class of monodromy representations
relative to different points in $B$ and different choices of identification
of $T_b B$ with $(\bR^n,\bZ^n)$.

The link between the topological and affine monodromies can be seen most
easily in local action-angle coordinates $(p,q)$ on a neighborhood of
a regular fiber $F_b=\pi^{-1}(b)$. 
The vectors $\frac{\del}{\del p_1}, \dots, \frac{\del}{\del p_n}$ at $b$
form a basis for $\Lambda_b$ and 
the homology classes of integral curves tangent to the vector fields
$(\frac{\del}{\del q_1},\ldots \frac{\del}{\del q_n})$ on $\pi^{-1}(b)$ 
represent a basis for $H_1(F_b,\bZ)$.
With respect to these bases, if the topological monodromy of the 
Lagrangian fibration along a loop $\gamma$ based at $b$
is given by $A\in GL(n,\bZ)$,
then the affine monodromy along $\gamma$ is
given by its inverse transpose $(A^{-1})^T$.
This follows immediately from the requirement that the endomorphism of 
$T_xM$, $x\in F_b$, determined by the topological and affine 
monodromies be symplectic.

\subsection{Neighborhood of a nodal fiber} \label{nodalnbhd.sec}

We now turn our attention to almost toric fibrations in dimension four,
the lowest dimension in which nodal fibers can occur.  (The reader
can assume for the rest of the paper that the dimension of the total
space is four.)

Since a neighborhood of a nodal fiber with one singular point is smoothly
equivalent to the fibered neighborhood of a singular fiber in a Lefschetz
fibration, the monodromy around the fiber is, with respect to some basis
for the first homology of a regular fiber $F_b=\pi^{-1}(b)$, 
\be
A_{(1,0)}:= \left(  
\begin{array}{cc} 
1 & 1 \\  
0 & 1 
\end{array} 
\right).
\ee
The reader should note that $A_{(1,0)}$ 
is a parabolic matrix with eigenvector 
$\left(\begin{smallmatrix}1\\0\end{smallmatrix}\right)$.
With respect to an arbitrary basis the monodromy matrix has the form
\be
A_{(a,b)}:= \left(  
\begin{array}{cc} 
1-ab & a^2 \\  
-b^2 & 1+ab 
\end{array} \right)
\ee
with eigenvector $\left(\begin{smallmatrix}a\\b\end{smallmatrix}\right)$
for some relatively prime $a,b\in \bZ$.
Viewing the singular fiber as a regular fiber with a circle pinched to
a point, this circle (which represents the {\it vanishing cycle}) 
represents the homology class $(a,b)$.

If there are $k$ singular points 
in a nodal fiber then the monodromy around such
a fiber is the product of $k$ nodal monodromy matrices all of which
have the same eigenvector (since the fiber is obtained from a regular fiber
by pinching $k$ circles all of which represent the same homology class).
Thus the affine monodromy matrix around a node of multiplicity $k$
is, with respect to some basis, 
\be
A^k_{(1,0)}= \left(  
\begin{array}{cc} 
1 & k \\  
0 & 1 
\end{array} 
\right).
\ee

While nodal fibers occur naturally in certain examples coming from algebraic
geometry and integrable systems, these typically do not give a clear
picture of the local fibered structure.  For an explicit local model
consult Section~4.4 of~\cite{Symington.fourfromtwo}.

It is important to note that while the germ of a 
neighborhood of a nodal fiber with a fixed number of
singular points is unique up to symplectomorphism,
it is not unique up to fiber-preserving symplectomorphism.
Indeed, Vu Ngoc, S.~\cite{VuNgoc.focusfocus}
has identified a 
non-trivial invariant that classifies the germs of such neighborhoods
up to fiberwise symplectomorphism.

Given our interest in the symplectic and topological properties
of the total space of an almost toric four-manifolds, the following uniqueness
statement suffices. 
\begin{prop} \label{nodalgerm.prop}
Consider symplectomorphic neighborhoods of a pair of nodal fibers in
an almost toric fibration.   The 
symplectomorphism between them can be chosen to be fiber-preserving on
the complement of smaller fibered neighborhoods.  
\end{prop}
The proof of this proposition given in~\cite{Symington.grbd}
can easily be modified to accommodate multiple singularities 
on the nodal fiber. 

\subsection{Almost toric bases in dimension four}

In dimension four an almost toric fibration can have three types of
singular points: elliptic of corank one or two, or nodal singular points.
The construction of toric manifolds given in Section~\ref{toric.sec}
and the normal form for elliptic singular points implies that the
image of any such point has a neighborhood
that is integral affine isomorphic to a neighborhood of the origin in
either the first quadrant of $\bR^2$ or the right half plane in $\bR^2$.
Meanwhile, the normal form for a nodal singularity
implies that the images of nodal fibers are isolated singularities of
the affine structure on the base.  The structure of this singularity is 
constrained
by the topological monodromy around the singular fiber which in turn
determines the affine monodromy around the node.

Specifically, suppose the monodromy around a nodal fiber, with respect to some
basis for the first homology of a regular fiber is
$A_{(a,b)}$.  Then the discussion of Section~\ref{monodromy.sec} implies that
the affine monodromy is $A_{(-b,a)}$.
Therefore the vector $(-b,a)$ is tangent to the one well-defined line
that passes through the node.
Accordingly we call this line the {\it eigenline} through the node.

Knowing the monodromy around an isolated singular point in an affine surface
does not completely determine the germ of its neighborhood.  
In particular, there is an infinite family of isolated singularities 
around which the monodromy is parabolic.  
To distinguish between them,
remove an {\it eigenray} $R$ based at the node from a neighborhood $N$ of
the node, choose a projection
$\Phi: (N-R,\cA)\ra (\bR^2,\cA_0)$ and count the number of preimages of
a generic point in the image.   
The following lemma, together with the fact that the monodromy around
a node is parabolic, implies the uniqueness of 
the germ of a neighborhood of a node. 

\begin{lemma} \label{nodenbhd.lem}
Let $N$ be a contractible neighborhood of a node in an affine two-manifold.
Suppose $N$ contains just one node and $R$ is a ray based at the node $b$ such
that $N-R$ is simply connected.   Then any immersion 
$\Phi: (N-R,\cA)\ra (\bR^2,\cA_0)$ is an embedding. 
\end{lemma}
A proof of this lemma can be found in~\cite{Symington.grbd} or Section 9.2
of~\cite{Symington.fourfromtwo}.  It relies on the fact (due to Gromov and
Eliashberg) that a fillable contact three-manifold is tight.
By the same reasoning as in the proof of Lemma~\ref{nodenbhd.lem}, the
germ of a neighborhood of a node of any multiplicity is completely determined
by the monodromy.

We now have that the neighborhood of any fiber in an almost toric
fibration, singular or not, has a neighborhood that can be recovered
(up to a variation in the fibration near a nodal fiber)
from the base of the fibration.  
Zung's study of Lagrangian fibrations with topologically stable 
non-degenerate singularities (\cite{Zung.II}) focused on how such 
neighborhoods can fit together.  An essential invariant is the
{\it Lagrangian Chern class}, an element of the first homology
of the base with values in the sheaf of closed basic one-forms
(one-forms that vanish on vectors tangent to fibers) modulo
those forms that arise from contracting the vector field for
a symplectic fiber-preserving circle action with the symplectic form.
(To be precise, this Chern class is actually a relative class in the sense 
that it is defined relative to a given reference fibration.)
From his work we extract the following
generalization of Proposition~\ref{toricbase.prop}:

\begin{prop} \label{basedetermination.thm}
In dimension four, an almost toric manifold is determined, 
up to fiber-preserving 
symplectomorphism, by its base $(B,\cA,\cS)$, the Lagrangian
Chern class and the local structure of the fibered neighborhoods of
its nodal fibers.  
\end{prop}
Note that if the base has 
the homotopy type of a zero or one-dimensional manifold, then
the Lagrangian Chern class vanishes.
Therefore Proposition~\ref{nodalgerm.prop} implies
\begin{corollary} \label{basedetermination.cor}
If an almost toric base $(B,\cA,\cS)$ of dimension two has
the homotopy type of a zero or one-dimensional manifold, then
$(B,\cA,\cS)$ determines up to symplectomorphism 
the total space of an almost toric fibration.
\end{corollary}

Throughout the rest of the paper we abuse language a bit and
refer to the almost toric manifold defined by a base.  More precisely,
the base diagram merely defines a symplectic manifold and a set
of almost toric fibrations that share the same base.

\subsection{Base diagrams and branch moves} \label{atbase.sec}

In this section we restrict ourselves to dimension four and
introduce base diagrams that 
allow for reconstruction of almost toric bases that immerse in $\bR^2$.
In light of Corollary~\ref{basedetermination.cor}, these can be viewed as
a generalization of moment maps.

\begin{defn} \label{branch.defn}
Consider an integral affine surface $(B,\cA)$ with nodes $\{b_i\}_{i=1}^k$
and non-empty boundary.
A set of {\it branch curves} for $(B,\cA)$ is a union of disjoint curves,
$R=\cup \{R_i\}_{i=1}^k$, 
such that each $R_i$ has one endpoint at $b_i$ and $R_i\cap \bdy \overline B$
is one point. 
\end{defn}
(Note that since $B$ is locally modeled on $\bR^2$ on the complement
of the nodes, $\overline B$ is well defined.)
An essential feature of branch curves is that whenever $B$ is a disk,
there is an immersion of $(B-R,\cA)$ into $(\bR^2,\cA_0)$.
Consequently, whenever the universal cover of
$B$ embeds smoothly in $\bR^2$, there exists an affine immersion 
into $(\bR^2,\cA_0)$ of a fundamental domain for $(B-R,\cA)$.

\begin{defn} \label{basediag.defn}
Suppose $(B,\cA,\cS)$ is 
an almost toric base with $B$ homeomorphic to a two-disk.
Let $R$ be a set of branch curves and
$\Phi:(B-R,\cA)\ra(\bR^2,\cA_0)$ an affine immersion.
A {\it base diagram} of $(B,\cA,\cS)$ with respect to $R$ and $\Phi$ is
the image $\Phi(B-R)$ with the following additional data:
\ben
\item Any portion of the boundary belonging to the closure of the 
one-stratum is drawn with heavy lines. 
\item An asterisk indicates the location of a node. 
\item Any information needed to recover the base from its image in $\bR^2$ if
it is not embedded.  
\een
If $B$ is not a two-disk but its universal cover does smoothly embed in 
$\bR^2$, then by a base diagram for $(B,\cA,\cS)$ we mean an affine image
of a fundamental domain for $B-R$ with the additional data as above.
\end{defn}
Since the third type of data is not usually necessary we do not set
a methodology for encoding it.  Often the geometry of the base is clear, 
even if it is not embedded, as in Example~\ref{exoticR4.ex}.

\begin{example} \label{nodenbhd1.ex}
Let $U$ be any neighborhood of the origin in $\bR^2$ and let $R$ be a 
ray with rational slope $\frac{b}{a}$ based at the origin.  (Assume
$a,b$ are relatively prime integers).
Then $U-R\subset\bR^2$ with an asterisk at the origin 
is a base diagram for a neighborhood of a nodal fiber
$\pi:(N,\omega)\ra (B,\cA)$.
Furthermore, with respect to coordinates (on the complement of a ray in
$(B,\cA)$) induced from the base diagram, the affine monodromy around the 
node is $A_{(a,b)}$.
\end{example}

The next example shows how a base diagram varies depending on 
the choice of ray that is removed.

\begin{example} \label{nodenbhd3.ex}
Suppose $\pi:(N,\omega)\ra (B,\cA)$ is an almost toric fibration of a 
neighborhood of a nodal fiber. 
Choose a ray $R$ based at the node and integral 
affine coordinates on the complement of $R$ such that the affine monodromy
is $A_{(1,0)}$.  If the ray $R$ belongs to a line with slope $\frac{b}{a}$
then there is a projection $\Phi:(N-R,\cA)\ra (\bR^2-S,\cA_0)$ where 
$S$ is the sector bounded by the vectors $(a,b)$ and $(a+b,b)$.  This
projection will be surjective onto $U-S$ 
where $U$ is a neighborhood of the origin.
\end{example}

Varying the choice of ray or curve that is removed in order to project
to $(\bR^2,\cA_0)$ constitutes a {\it branch move}.  Thus there are
two ways to vary a base diagram when $B$ is a two-disk: 
via branch moves and by changing the
projection by composing with an element of ${\rm Aff}(2,\bZ)$.  (Of course,
if $B$ is not a two-disk, then we can also vary the choice of fundamental
domain.)
While one base that immerses into $\bR^2$ has many base diagrams, from
any one of them one can reconstruct the base.

\section{Almost toric surgeries}

\subsection{Nodal trades} \label{nodaltrade.sec}

In this section we describe a surgery operation that changes an almost
toric fibration of a symplectic four-manifold into another almost toric 
fibration of
the same symplectic manifold.  The essential idea is that one can trade
a zero-dimensional singular fiber for a nodal fiber (and vice versa under
the appropriate conditions).

We start with an even simpler way to modify almost toric fibrations:

\begin{defn} \rm \label{nodalslide.def}
Two almost toric bases $(B,\cA_i,\cS_i)$, $i=1,2$, are related by
a {\it nodal slide} if there is a curve $\gamma\subset B$
such that $(B-\gamma,\cA_1,\cS_1)$ and $(B-\gamma,\cA_2,\cS_2)$ 
are isomorphic and, for each $i$, $\gamma$ contains one node of
$(B,\cA_i,\cS_i)$ and $\gamma$ belongs
to the eigenline through that node.
\end{defn}

A nodal slide should be thought of as a one-parameter family of
almost toric bases in which a node moves in the base along its eigenline.  
Of course it corresponds to a one-parameter family of almost toric fibrations 
of one manifold.  Exactness of the symplectic structures on the preimage of
a neighborhood of $\gamma$  allows us to use Moser's argument
to confirm that the symplectic manifolds that fiber over
$(B,\cA_i,\cS_i)$, $i=1,2$, are symplectomorphic.  This perspective
on nodal slides makes it easy to find a one-parameter family of almost
toric fibrations connecting a fibration with a singular fiber having
$k$ nodes to a fibration in which each fiber has only one singular point.  
In the base, at one extreme one would have 
a node of multiplicity $k$ and at the other
one would have $k$ nodes that live on one line, the eigenline.

If the eigenline through a node intersects the one-stratum of the
base then the limit of nodal slides as one endpoint approaches this
stratum will
result either in changing the topology of the total space to form
an orbifold, or else merely a change in the fibration that replaces the
nodal fiber with an elliptic singular point of corank two.  We call this
operation a {\it nodal trade}.  Zung~\cite{Zung.II} had observed that 
this operation could be performed on Lagrangian fibrations and 
that the one-parameter family connecting the initial and final fibrations 
appears frequently in integrable systems: it is a Hamiltonian-Hopf 
bifurcation~\cite{vanderMeer.HamiltHopf}.

\begin{figure}
\begin{center}
   	\psfragscanon
	\psfrag{a}{(a)}
	\psfrag{b}{(b)}
	\includegraphics[width=4.5in]{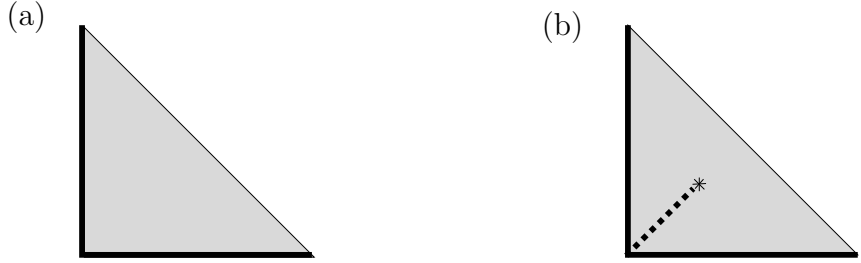}
\end{center}
\caption{Symplectomorphic four-balls.}
\label{balls.fig}
\end{figure}

As an example of a pair of bases related by a nodal trade, consider
the base diagrams shown in Figures~\ref{balls.fig} (a) and (b).
These base diagrams define symplectomorphic manifolds.

Indeed, Figure~\ref{balls.fig} (a) is a base diagram corresponding to a toric
fibration of the standard symplectic four-ball.
Meanwhile, Figure~\ref{balls.fig} (b) 
is a base diagram for an almost toric four-ball: 
the preimage of a collar neighborhood
of the one-stratum (not including the node) 
is an $S^1\times D^3$, or $D^4$ with a one-handle attached; the preimage
of the whole base differs by attaching a $-1$-framed two-handle that
is a thickening of the vanishing disk of the nodal fiber, but this two-handle
cancels the one-handle yielding $D^4$.

Since the base diagram in Figure~\ref{balls.fig} (a) is the limit
of a nodal slides of the node in  Figure~\ref{balls.fig} (b), we can
again invoke Moser's argument to establish that the total spaces fibering
over the two bases are symplectomorphic.

Here is a precise definition:
\begin{defn}
Two almost toric bases $(B_i,\cA_i,\cS_i)$, $i=1,2$, 
differ by a {\it nodal trade}
if each contains an arc $\gamma_i$
such that $(B_i-\gamma_i,\cA_i,\cS_i)$, $i=1,2$, 
are isomorphic, and the cardinality
of the zero stratum of $(B_1,\cA_1,\cS_1)$ is one less than 
that of $(B_2,\cA_2,\cS_2)$.
\end{defn}

The argument that the base diagrams in Figure~\ref{balls.fig} define
symplectomorphic manifolds generalizes in the obvious way to hold 
for any nodal trade.  
Therefore,

\begin{thm} \label{nodaltrade.thm}
Two almost toric bases that are related by a nodal trade are symplectomorphic.
\end{thm}

\subsection{Rational blowdowns and generalizations} \label{blowdown.sec}

The rational blowdown is a surgery in which the neighborhood of a chain
of spheres with lens space boundary $L(n^2,n-1)$, or more 
generally $L(n^2,nm-1)$,
is replaced by a manifold with the same boundary but having 
rational homology equal to that of a four-ball.
This surgery, useful in the study of smooth four-manifolds, was introduced
by Fintushel and Stern~\cite{FintushelStern.blowdown} and its generalization
by Park~\cite{Park.gblowdowns}.
The second author proved that these surgeries can be performed in the 
symplectic category~(\cite{Symington.blowdowns},\cite{Symington.grbd}),
thereby showing certain exotic four-manifolds could be symplectic.
The proof relied on the fact that the collar neighborhoods
involved in the surgery admit a toric fibration.

\begin{obs} The generalized rational blowdown is an almost toric surgery:
the manifolds removed and glued in are both almost toric (in fact the
former, a neighborhood of a chain of symplectic spheres, can always 
be chosen so that it is toric).
See Example~\ref{rationalballs.ex} for the almost toric structure on the
rational ball that gets glued in. 
In this setting one can always assume that the gluing locus is a neighborhood
of a contact manifold (or equivalently, that the boundaries in question are
contact).
\end{obs} 

This almost toric perspective leads to a generalization in which 
one exchanges almost toric manifolds whose boundaries may not be contact,
but nonetheless have collar neighborhoods that are toric.
To be useful as a surgery, 
one would need a method for finding embedded lens spaces that
bound and have a toric neighborhoods.  
Generalizing the phenomenon that any chain of symplectic surfaces has
a toric neighborhood, one could hope for the following:

\begin{conj}
Suppose a symplectic four-manifold $(M,\omega)$ contains a top-dimensional
submanifold with boundary that is diffeomorphic to an almost toric manifold 
whose base has a one-stratum diffeomorphic to an interval.  Then $(M,\omega)$
contains an almost toric submanifold with boundary diffeomorphic to and
contained in the given one. 
\end{conj}

Sets of almost toric manifolds that could then be exchanged are then provided
by: 

\begin{prop} \label{surgery.prop}
Consider two sequences of nodal monodromy matrices, \\
$\{A_1,\ldots, A_m\}$ and $\{B_1,\ldots, B_n\}$.
Suppose there is a vector $v$ such that 
\be \label{match.eq}
A_mA_{m-1}\cdots A_1 v = B_n B_{n-1} \cdots B_1 v.
\ee
Suppose also that the number of times
$A_i A_{i-1} \cdots A_1 v \times v$ changes sign as $i$ ranges from $1$ to $m$
equals the number of times $B_j B_{j-1}\cdots B_1v\times v$ changes
sign as $j$ ranges from $1$ to $n$.
Then there are two almost toric manifolds $(M,\omega_M)$ and $(N,\omega_N)$ 
such that
\ben
\item their bases have empty zero-strata,
\item their boundaries have collar neighborhoods that are 
fiber-preserving symplectomorphic, and
\item the matrices 
$A_i$ and $B_j$ are the monodromy matrices across branch curves in
base diagrams for each.
\een
\end{prop}

\begin{proof}
Let $U_M$ be the one-sided closed immersed neighborhood in $\bR^2$
of a polygonal line $P_M$
such that the tangent vectors to the (ordered) linear components are
$v$, $A_1v$, \ldots $(A_mA_{m-1}\cdots A_1) v$ and the boundary of $U_M$ is
locally convex on this polygonal part.  (See Figure~\ref{U_M.fig} in
which $P_M$ is drawn with a heavy line.)

\begin{figure}
\begin{center}
        \psfrag{v}{$v$}
        \psfrag{A1v}{$A_1v$}
        \psfrag{A2A1v}{$A_2A_1v$}
        \psfrag{A7A2A1v}{$A_7\cdots A_2A_1v$}
	\includegraphics[width=2.5in]{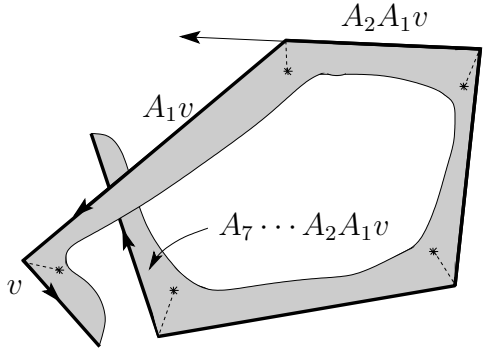}
\end{center}
\caption{Neighborhood $U_M$.}
\label{U_M.fig}
\end{figure}

Indicate $m$ nodes that are the endpoints of disjoint line segments,
each with one endpoint at a vertex of the polygonal boundary and lying
in an eigenline for the corresponding monodromy matrix (that relates
the tangent vectors of the edges that meet at the given vertex).
Do the same for the sequence of vectors 
$v$, $B_1v$, \ldots $(B_nB_{n-1}\cdots B_1)v$, creating the immersed
domain $U_N$ that is a one-sided neighborhood of a polygonal line $P_N$.
Equation~\ref{match.eq} guarantees that the boundaries of the almost toric
manifolds defined by these two base diagrams are diffeomorphic.
The additional fact that the polygonal lines $P_M,P_N$ 
wind around the origin the same number
of times then assures that one can choose the lengths of the edges in 
the polygonal boundaries so that $P_M$ and $P_N$ differ by an element
of ${\rm Aff}(2,\bZ)$.  Consequently, we can assure that the collar 
neighborhoods
of the boundaries of the manifolds fibering over $U_M,U_N$ are 
symplectomorphic.
We are of course appealing to Corollary~\ref{basedetermination.cor} in passing
from properties of the base diagrams to properties of the total spaces. 
\end{proof}
Whenever $m\ne n$ the manifolds $M,N$ have different Euler characteristics so
performing the surgery would definitely cause a change in topology.

It would be interesting to see how the class of almost toric manifolds
with lens space boundary compares with Lisca's complete 
list~\cite{Lisca.lensfillings} of symplectic
fillings of lens spaces with standard contact structures (coming from the
standard tight structure on $S^3$ via a finite group action).

\subsection{Almost toric blow-ups} \label{atblowup.sec}

Any locally toric manifold can be blown up at a point, yielding a
locally toric fibration of the blow-up, provided the base 
has non-empty zero-stratum.
In the almost toric category it suffices for the one-stratum
to be nonempty.  (This fact was first observed by Zung~\cite{Zung.II} in
the context of four-dimensional singular Lagrangian fibrations.)

\begin{figure}
\begin{center}
   	\psfragscanon
	\psfrag{a}{(a)}
	\psfrag{b}{(b)}
	\includegraphics[width=4.5in]{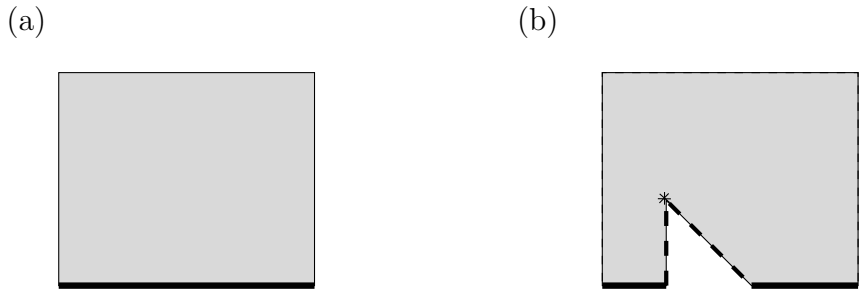}
\end{center}
\caption{$S^1\times D^3$ and $(S^1\times D^3)\,\#\,\overline{\bC P}^2$. }
\label{blowup.fig}
\end{figure}

Indeed, consider the base diagrams shown in Figures~\ref{blowup.fig} (a) and
(b).  The first defines a toric fibration of $S^1\times D^3$ and 
the second represents the base of an almost toric fibration of 
$(S^1\times D^3)\,\#\,\overline{\bC P}^2$.  The $-1$-sphere that is introduced
in the blow-up can be found in the preimage of an arc
connecting the node and the one-stratum.  See Section 5.4
of~\cite{Symington.fourfromtwo} for more details.

Theorem~\ref{basedetermination.thm} then 
implies that up to scaling, the action of
${\rm Aff}(2,\bZ)$, or a branch move, Figures~\ref{blowup.fig} (a) and (b)  
completely describe an almost toric surgery that amounts to blowing up
the total space.
Generalizing in an obvious way, we define:

\begin{defn}
An {\it almost toric blow-up} of an almost toric fibered manifold
$\pi:(M,\omega)\ra(B,\cA,\cS)$ is an almost toric fibration
$\pi':(M\,\#\,\overline{\bC P}^2,\omega')\ra (B',\cA',\cS')$ such that there is
an arc $\gamma'\subset B'$ based at a node and such that 
$(B'-\gamma',\cA',\cS')$  embeds in $(B,\cA,\cS)$.
\end{defn}

\section{Almost toric manifolds: examples}\label{examples.sec}

\subsection{Non-compact/non-empty boundary} \label{noncpct.sec}

We start with an example which is toric but which has a base that
does not (affinely) embed in $(\bR^2,\cA_0)$.  
This example was first noticed by Zung~\cite{Zung.II} and was inspired
by an example due to Bates and Peschke~\cite{BatesPeschke.exotic}.

\begin{figure}
\begin{center}
	\includegraphics{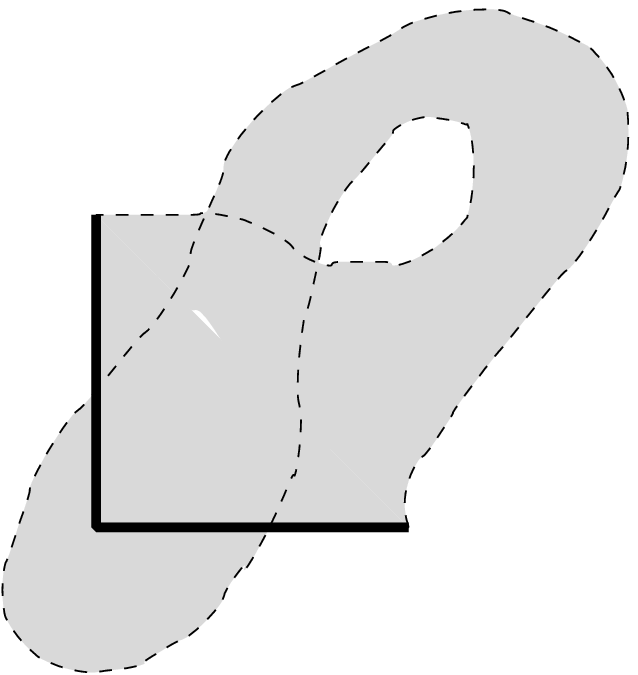}
\end{center}
\caption{Exotic $\bR^4$.}
\label{exoticR4.fig}
\end{figure}

\begin{example} \rm \label{exoticR4.ex}
Figure~\ref{exoticR4.fig} is a base diagram for a toric fibration
of an open symplectic four-ball that is exotic: there is no symplectic
embedding of this ball into $(\bR^4,\omega_0)$ where 
$\omega_0$ is the standard symplectic structure.
In this figure the vertex (the zero-stratum of the base) is at the origin.

The total space can easily be seen to be diffeomorphic to an open ball 
since the base has the same topology and stratification as the base diagram
corresponding to the standard moment map image
of a neighborhood of the origin in $(\bR^4,\omega_0)$.  Indeed,
the base diagram of the exotic four-ball differs 
only by having a lobe stretched out (in $\bR^2$) and wrapped around 
until it covers the origin.  Throughout, the total space fibering over
the lobe is diffeomorphic to $T^2\times D^2$.

To see that the symplectic structure $\omega$ is exotic we appeal to
Gromov's theorem that there are no exact
Lagrangian tori in $(\bR^{2n},\omega_0)$~\cite{Gromov.pseudohol}.
Having placed the vertex at the origin, we 
have that $\alpha=p\,dq$ 
(with respect to the 
induced local coordinates on the preimage of the top dimensional stratum) 
is a primitive of the symplectic structure which 
extends to the preimage of the lower-dimensional strata, i.e. is a global
primitive.
The preimage of the origin is the union of a point and an
exact Lagrangian torus $T$ (since $\alpha|_T=0$ and hence represents
a trivial class in cohomology).
\end{example}

\begin{example} \label{pendulum.ex}
The phase space of a spherical pendulum in a gravitational field
is an integrable system that has one nodal fiber.   
Indeed, the energy and angular momentum functions form a map to $\bR^2$
that defines a Lagrangian fibration.
The singular fibers consists of:
\ben
\item A point corresponding to when the pendulum is at its lowest
position with zero kinetic energy.
\item A nodal fiber, the singular point of which corresponds to the
pendulum at its highest position with zero velocity; the rest of the
nodal fiber is made up of orbits in which the pendulum swings in a
vertical plane with the same total energy as the singular point.
\item A one parameter family of circle fibers each of which is an orbit
in which the pendulum rotates around the central vertical axis with maximal
angular momentum for a given energy.
\een
Note that using the energy and angular momentum as the components of the
fibration map, the induced affine structure on $\bR^2$ is not the standard
one.

For more details the reader can consult 
Duistermaat~\cite{Duistermaat.actionangle}.
\end{example}

\begin{example} \label{rationalballs.ex}
A rational ball (four-manifold with boundary having the same rational homology
as a four-ball) admits a simple almost toric fibration whose singular fibers
consist of a one-parameter family of circle fibers and a nodal fiber.  They
have base diagrams that look very similar to the moment map images of 
symplectizations of lens spaces with an $S^1$ invariant contact structure. 
In fact, the only difference is the presence of a
dashed line indicating a branch curve, and hence a node. 
See~\cite{Symington.fourfromtwo} for more details.

With this perspective it is easy to check that portion of
the phase space of the spherical pendulum with bounded the energy (for a
sufficiently large bound) is an almost toric rational ball.
with boundary $L(2,1)$. 

Note that the four-ball itself admits the same type of fibration as a
rational ball: start with the standard toric fibration and perform a
nodal trade.

Interest in rational balls stems from  
their role in constructing four-manifolds via 
the rational blowdown surgery described briefly in Section~\ref{blowdown.sec}.
\end{example}

\subsection{Closed manifolds} \label{closed.sec}

In this section we describe almost toric fibrations of the following 
manifolds:
$S^2\times S^2$, $\bC P^2\,\#\,n\overline{\bC P}^2$, 
$(S^2\times T^2)\,\#\, n\overline{\bC P}^2$, $S^2\tilde\,\#\, T^2$, the K3 surface,
the Enriques surface, and
certain torus bundles over $T^2$ and the Klein bottle.   
Note that this is the full set of manifolds that appears in 
Table~\ref{diffclass.table}.

A majority of these manifolds admit complex structures.  In the
case of toric fibrations of the the rational surfaces and 
almost toric fibrations of the K3 surface and $T^4$ the fibrations arise
naturally in algebraic geometry: 
the rational surfaces are precisely the
toric algebraic surfaces while the K3 surface and $T^4$ are examples
of hyperk\"ahler manifolds that, via a hyperk\"ahler rotation admit both
holomorphic and Lagrangian fibrations.   
The torus bundles over tori
that admit Lagrangian fibrations include the Kodaira manifolds
that were the first examples of non-K\"ahler symplectic manifolds
(Thurston~\cite{Thurston.nonKahler}).
Note that while some of the torus bundles over tori have $b_1=2$, these
manifolds do not admit a complex structure; in particular, the
hyperelliptic surfaces do not admit regular Lagrangian fibrations.
At the end of this section we prove that an almost toric manifold
fibering over $\bR P^2$ is diffeomorphic to an Enriques surface.
Because the fibration is not compatible with a complex structure, we do this 
by showing that
the Lagrangian-fibered manifold can be obtained from an almost toric 
fibration of the elliptic surface $E(1)$ 
by performing two smooth log transforms of multiplicity two.

\begin{example} \label{rational.ex}
The rational surfaces $\bC P^2\,\#\,n\overline{\bC P}^2$ and 
$\bC P^1\times\bC P^1$ 
are all of the symplectic four-manifolds that admit a toric structure.
In addition to toric fibrations these manifolds admit almost toric fibrations
that contain nodal fibers and hence cannot be toric.  Indeed, one can
perform nodal trades (Section~\ref{nodaltrade.sec}) to replace
zero-dimensional fibers with nodal fibers.  

One nice feature of the almost toric fibrations is that they provide 
families of fibrations that interpolate between toric fibrations of the same
symplectic manifold.  (See~\cite{Symington.fourfromtwo} Section 6.2.)  
This is always true
for Hirzebruch surfaces (diffeomorphic to $S^2\times S^2$ or 
$\bC P^2 \,\#\, \overline{\bC P}^2$) and we conjecture that it is true for
any symplectic manifold that admits a toric structure. 
\end{example}

Hyperk\"ahler surfaces are equipped with a two-sphere of complex
and corresponding K\"ahler structures.  This allows, via a hyperk\"ahler
rotation, to transform a holomorphic fibration into a special
Lagrangian fibration 
(a Lagrangian fibration that is also adapted to the complex 
structure).  
In complex dimension two the only closed hyperk\"ahler manifolds are
the K3 surfaces and complex tori ($T^4$ equipped with a complex structure).  

\begin{example} \label{K3.ex}
As a complex manifold there are many elliptically fibered K3 surfaces,
and generically they have only nodal singular fibers.
Performing a hyperk\"ahler rotation on such a generic 
K3 surface yields an almost toric fibered K3 with $24$ nodal fibers.

Note that all K3 surfaces are diffeomorphic so among smooth
four-manifolds there is just one K3 surface.  
One way to construct an almost toric fibration of the K3 surface is
inspired by the well known fact that the K3 surface, or $E(2)$, is
the fiber sum of two copies of the elliptic surface
$E(1)=\bC P^2\,\#\, 9\overline{\bC P}^2$ which fibers over $\bC P^1$ with
$12$ nodal fibers.  Specifically,
\ben
\item Choose two toric fibered copies of $\bC P^2\,\#\, 9\overline{\bC P}^2$.
\item The base diagram of each has twelve vertices.  Perform nodal
trades at all of the vertices.  This yields two almost toric manifolds,
each with
a smooth symplectic torus of self-intersection zero fibering over the
boundary of the base.  This almost toric
fibration has $12$ singular fibers that are 
Lagrangian and hence are not the singular fibers of a symplectic or
holomorphic fibration.
However, the preimage of the boundary is a symplectic torus that can be
viewed as a regular fiber of a symplectic fibration.
\item Symplectic sum the two almost toric copies of
$\bC P^2\,\#\, 9\overline{\bC P}^2$
along these symplectic tori.  In the
base this amounts to joining the two bases along their  boundaries.  
\een 
\end{example} 

\begin{example} \label{toriovertori.ex}
Regular Lagrangian fibrations over $T^2$ have been classified up to 
fiber preserving symplectomorphism by
Mishachev~\cite{Mishachev.Lagbundles}.  
Since an almost toric fibration over a torus cannot 
have singular fibers (Theorem~\ref{base.thm}), there are no other examples
fibering over the torus.
 
For identifying the possible total spaces, it is most convenient
to turn to Geiges' classification of cohomology classes with symplectic
representatives on torus bundles
over tori~\cite{Geiges.torusbundles}.
In that paper he specifies which torus bundles 
admit symplectic or Lagrangian fibers. 

Recall that up to fiber-preserving diffeomorphism a torus bundle over a torus
is specified by two monodromy matrices and two integers that are
the obstruction to the existence of a section and can be viewed as 
a Chern class.

From~\cite{Geiges.torusbundles} we extract the following:
\begin{lemma} \label{Geiges.lem}
If a torus bundle over a torus admits a Lagrangian fibration then it
has monodromy 
\be
\left\{I,\left(  
\begin{array}{cc} 
1 & \lambda \\  
0 & 1 
\end{array} 
\right)\right\} \quad {\rm with\ } \lambda\in\bZ
\ee
and arbitrary Chern class $(m,n)\in\bZ^2$.
\end{lemma}

Referring then to~\cite{SakamotoFukuhara.torusbundles} we see that the torus
bundles specified by the integers $(\lambda,m,n)$ and $(\lambda',m',n')$ are
equivalent if and only if $\lambda'=\epsilon\lambda$ and $n'=\epsilon n$ where
$\epsilon\in\{1,-1\}$,  and $m'= m +k\lambda +l n$ for some integers $k,l$.
Furthermore, if $\lambda=0$ or $n=0$ (or equivalently, if $b_1\ge 3$),
then the total space is diffeomorphic
to the total space of a fibration determined by $(\mu,0,0)$ for some
$\mu$; otherwise the diffeomorphism classification agrees with 
the bundle classification.
\end{example}

\begin{figure}
\begin{center}
   	\psfragscanon
	\psfrag{a}{(a)}
	\psfrag{b}{(b)}
	\includegraphics[width=4.5in]{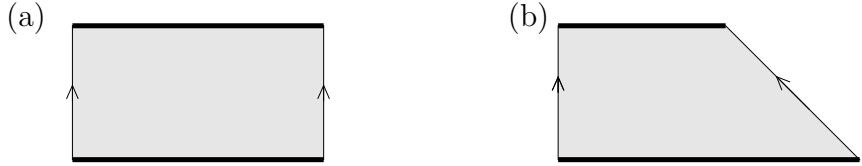}
\end{center}
\caption{Sphere bundles over tori.}
\label{sphereovertorus.fig}
\end{figure}

\begin{example} \label{sphereovertorus.ex}
Up to diffeomorphism there are two sphere bundles over the torus.  Both
of these admit an almost toric fibration over the cylinder. 
Base diagrams of their fundamental domains are shown in 
Figure~\ref{sphereovertorus.fig}.
The first one clearly corresponds to $S^2\times T^2$.  (Note that the preimage
of any vertical cross-section of the base is a copy of $S^2\times S^1$.)
To see that the second one corresponds to a non-trivial bundle, one can
detect from the base diagram that the preimage of the lower component
of the one-stratum is a torus of self-intersection $1$.  Details 
can be found in Section $7$ of~\cite{Symington.fourfromtwo}.

Blow-ups of these manifolds also admit almost toric fibrations as
we can perform an arbitrary number of (sufficiently small)
almost-toric blow-ups.  Similar to sphere bundles over the sphere, 
blow-ups of the trivial and non-trivial sphere bundles over the torus
are diffeomorphic.  This can easily be proven in terms of base
diagrams by making a branch move, effectively switching the almost toric
blowup point from the torus in the preimage of one component of the one-stratum
to the torus in the preimage of the other component.
\end{example}
 
The base $B$ of an almost toric fibration, like a Lagrangian submanifold, need
not be orientable.  If the base is non-orientable then its double cover
$\tilde B$ is an integral
affine manifold with nodes that is again an almost toric 
base (since it is locally isomorphic to $B$).

\begin{example} \label{torioverKlein.ex}
If the base is a Klein bottle then the monodromy must be 
$\left\{\left( 
\begin{smallmatrix}
1 & 0 \\  
0 & -1 
\end{smallmatrix} 
\right), \left(\begin{smallmatrix}
1 & \lambda \\  
0 & 1 
\end{smallmatrix} 
\right)\right\}, \ \ \lambda\in\bZ.$ 
As in the case of a torus base, we can twist a fibration with section
by pulling out a fiber and regluing by a diffeomorphism such that the
section no longer extends through that neighborhood.  By lifting such
a twist to the double cover, one sees that the covering manifold must
be a torus bundle over a torus given by $(\lambda, m, n)$ where $m$ is
even and $n=0$.  In particular, this shows that only those torus
bundles over tori that have $b_1\ge 3$ admit a $\bZ_2$ action whose
quotient is a Lagrangian fibration over the Klein bottle. 
\end{example}

\begin{example} \label{torioverMoeb.ex}
To get an almost toric fibration over a M\"oebius band we need to 
take a fiber-preserving $\bZ_2$ quotient of an almost toric fibration over 
a cylinder.  The result will again be a sphere bundle over a torus.

Indeed, as when the base is a cylinder, a base diagram for a fundamental
domain
of a M\"oebius band base determines $S^2\times S^1\times I$ together with
some symplectic structure and an identification of the boundary components 
$S^2\times S^1\times \{0\}$ and 
$S^2\times S^1\times \{1\}$.
This identification differs from the cylinder case by precomposition with
a rotation by $\pi$ of the $S^2$ factor.  Since this map is isotopic to the
identification map in the cylinder case, the resulting total space is the 
same.
\end{example}
 
\begin{example} \label{Enriques.ex}
Any almost toric fibration over $\bR P^2$ is double covered by an almost
toric fibration over $S^2$, namely a K3 surface.  
This would lead one to suspect that such a manifold must be diffeomorphic
to the Enriques surface, a complex manifold that is a $\bZ_2$ quotient of
a K3 surface.  However, a holomorphic fibration of the Enriques surface
fibers over $\bC P^1$.
Since the $\bZ_2$ action in question is not holomorphic, 
we appeal to an argument in the
smooth category that is compatible with almost toric structures.
\end{example}

\begin{lemma} \label{Enriques.lem}
The total space of a Lagrangian fibration over $\bR P^2$ is a diffeomorphic to
an Enriques surface.
\end{lemma}

\begin{proof}
It is well known that the Enriques surface can be obtained from 
$E(1)=\bC P^2\,\#\,9\overline{\bC P}^2$ by performing two smooth multiplicity
two log transforms along fibers of an elliptic (or genus one
Lefschetz) fibration of $E(1)$ (cf.~\cite{GompfStipsicz.4mflds}).
The log transform, as a smooth operation, amounts to removing a 
neighborhood of a regular torus fiber 
and gluing it back in via a diffeomorphism of the boundary.
When performed on an elliptic fibration (or more
generally in the neighborhood of a cusp fiber) the effect of this operation 
on the total space depends only on an integer, the multiplicity: it
produces a multiple fiber of multiplicity $p$, namely
a fiber $f_p$ whose homology class satisfies
$p[f_p]=[f]$ where $f$ is a regular fiber.

The elliptic surface $E(1)$ has an almost toric fibration over the disk
with $12$ nodal fibers.  
In this fibration the boundary of the base $B$ has no zero-stratum, so the
preimage of the boundary is a smooth symplectic torus $T$.  This
torus can be viewed as a regular fiber of a Lefschetz fibration
of $E(1)$.
To get a fibration over $\bR P^2$ we can remove a 
neighborhood of $\bdy B$ that is fibered by geodesics 
parallel to the boundary (with respect to the affine structure)
and replace it with a M\"oebius band also fibered by geodesics parallel 
to the boundary.
This surgery on the base corresponds to a surgery of the total space
in which we remove a neighborhood of $T$ and glue
back in an almost toric fibration over a M\"oebius band, i.e.
the product of $S^1$ and the circle bundle over a M\"oebius band whose total
space is orientable.
This latter three-manifold can be obtained from a solid torus
by performing two Dehn surgeries of multiplicity two
along circles parallel to the core (cf.~\cite{Hatcher.threemanifolds}).
Since the the product of the identity map (on $S^1$) and a Dehn surgery
(on a solid torus) is a log transform we have that the surgery on the
base corresponds to performing two log transforms of
multiplicity two.  Since the product of the $S^1$ factor and circles
parallel to the core of the solid torus correspond to fiber tori of a
Lefschetz fibration of $E(1)$ we have that the resulting manifold is
indeed the Enriques surface.
\end{proof}

\section{Classification of closed almost toric manifolds}

\subsection{Possible bases} \label{possiblebases.sec}

Our goal in this section is to determine what affine surfaces with
nodes and stratification, i.e. what triples $(B,\cA,\cS)$, 
can be the base of a closed almost toric manifold.
Theorem~\ref{base.thm} establishes that the bases are precisely those
which appear in Table~\ref{diffclass.table} with the given number of
nodes and vertices, and that if the base is
a cylinder or M\"oebius band 
the eigenlines of any nodes must be parallel to the boundary.

The first step in proving Theorem~\ref{base.thm} is to determine that the
base $B$ must have non-negative Euler characteristic 
(Lemma~\ref{Eulerchar.lem}).
The essence of the argument is as follows: the integral
affine structure provides a flat structure on the complement of the nodes,
with respect to which
the boundary (if nonempty) is piecewise linear and locally convex;
then, in a rough sense, the nodes contribute non-negative curvature.
The curvature contributions cannot be measured using a metric compatible
with the affine structure as the presence of 
nontrivial affine monodromy is an obstruction to the existence of such
a metric.
However, on disks with nodes we can construct metrics that are inspired
by base diagrams, so called boundary compatible metrics
(Definition~\ref{bdycompat.def}).
These metrics allow us to bound the total turning angle as the boundary
is traversed counter-clockwise. 
The Gauss-Bonnet theorem thereby provides a lower bound on the total curvature
of the disk (Lemma~\ref{disk.lem}).  To rule out higher genus surfaces 
occurring as a base, the disk we work with is a fundamental domain in 
the universal cover.

\begin{defn} \label{bdycompat.def}
A metric $g$ on an integral affine disk with nodes $(D,\cA)$ is 
{\it boundary compatible} (with $\cA$) if a collar neighborhood of 
the boundary can be covered by a pair of open sets 
$\{U,V\}$ such that
\ben
\item $U$ is nonempty and simply connected,
\item $g|_{U}=\Phi^* g_0$ for some integral
affine map $\Phi:(U,\cA)\ra(\bR^2,\cA_0)$ and
\item $\bdy D\cap V$ is geodesic with respect to both $g$ and $\cA$
whenever $V$ is nonempty.
\een
\end{defn}

\begin{lemma} \label{bdycompat.lem}
If an integral affine disk with nodes $(D,\cA)$ has a boundary
that contains a linear segment (a subset homeomorphic to an interval and
geodesic with respect to $\cA$) then it admits a boundary compatible metric.
\end{lemma}

\begin{proof}
Let $\{b_i\}_{i=1}^k$ be the nodes.  Choose a set of branch curves $R$ whose
endpoints on the boundary all belong to the interior of one linear segment.
Let $V$ be a connected open subset that contains all these endpoints and
such that $V\cap\bdy D$ is a linear segment.
Let $U$ be an open simply connected set that covers $\bdy D-V$.
The simple connectedness of $U$ guarantees the existence of 
an integral affine immersion $\Phi:(U,\cA)\ra(\bR^2,\cA_0)$.
Let $g_U=\Phi^*g_0$ and let $g_V$ be a metric defined on $V$ such that
$\bdy D\cap V$ is geodesic with respect to $g_V$.  
Then construct a metric $g$ on $U\cup V$ 
from $g_U$ and $g_V$ using a partition of unity subordinate to 
$\{U,V\}$.  Extending $g$ arbitrarily to the rest of the disk, one obtains
a boundary compatible metric.
\end{proof}

If the boundary of $(D,\cA)$ contains no linear segments then
Definition~\ref{bdycompat.def} implies $U$ covers $D$, in which case there
could not be any nodes in $(D,\cA)$.
Requiring the existence of a linear segment when the monodromy is non-trivial 
facilitates comparison of the geodesic curvature along
$\bdy D$ with the geodesic curvature along the boundary of 
a flat disk that is isomorphic to the closure of $D-R$ (as in the
proof of Lemma~\ref{disk.lem}).
Note that in all our applications of boundary compatible metrics
there is no loss of generality to assume that $\bdy D$
contains a linear segment.

\begin{defn} \label{flat.def}
Given an integral affine disk with nodes $(D,\cA)$
consider an immersion $\Phi:(D-R,\cA)\ra(\bR^2,\cA_0)$ for
some set of branch curves $R$.
A {\it flat model} for $(D,\cA)$ with respect to $R$ and $\Phi$ is the 
flat disk
$(\widehat D,\hat g)$ where $\widehat D=\overline{D-R}$ and 
$\hat g=\Phi^* g_0$.
\end{defn}

\begin{lemma} \label{disk.lem}
Let $(D,\cA)$ be an integral affine disk with nodes equipped with a 
boundary compatible metric $g$.  Then $\int_D K_g \,dA\ge 0$. 
\end{lemma}

\begin{proof}
Since $g$ is boundary compatible, there is a pair of sets $\{U,V\}$ 
as in Definition~\ref{bdycompat.def} with $g|_U=\Phi^*g_0$.
Because $U$ is simply connected we can find a 
set of branch curves $R=\{R_i\}_{i=1}^k$ emanating
from the nodes $\{b_i\}_{i=1}^k$ that are disjoint from $U$, and
hence have their other endpoints on $\bdy D\cap V$.  (If there are no
nodes, then $R$ is empty.)  
Replace the affine immersion $\Phi$, defined on $U$, by its extension
to all of $D-R$.
Let $(\widehat D,\hat g)$ be the flat model for $(D,\cA)$ with respect to 
$R$ and $\Phi$, so $\hat g=\Phi^*g_0$ and $\hat g|_U=g|_U$.
The Gauss-Bonnet Theorem implies
\be
\int_D K_g\, dA = 2\pi-\int_{\bdy D} \kappa_{g}\, ds-\beta_{g} 
{\rm \ \ and\ \ }
0=2\pi-\int_{\bdy \widehat D} \kappa_{\hat g}\, ds-\beta_{\hat g}.
\ee
where $\kappa_g,\kappa_{\hat g}$ 
are the geodesic curvatures along the smooth parts of 
$\bdy D,\bdy\widehat D$ and $\beta_g,\beta_{\hat g}$ are the sums of the
turning angles at the vertices of $\bdy D,\bdy\widehat D$.  
By construction we have 
\be
\int_{\bdy \widehat D} \kappa_{\hat g}\, ds + \beta_{\hat g}=
\int_{\bdy D} \kappa_{g}\, ds + \beta_{g}+\sum_{i=1}^k \theta_i
\ee
where each $\theta_i$ is the contribution to the total turning angle along
the portion of $\bdy\widehat D$ introduced by the node $b_i$. 
For instance, Figure~\ref{turning.fig} shows a base diagram with a branch
cut emanating from a node with affine monodromy $A=A_{(1,-1)}$.  In this
case the node would contribute $\theta=\pi/4$ to the total turning angle.

\begin{figure}
\begin{center}
   	\psfragscanon
	\psfrag{v}{v}
	\psfrag{Av}{Av}
	\includegraphics{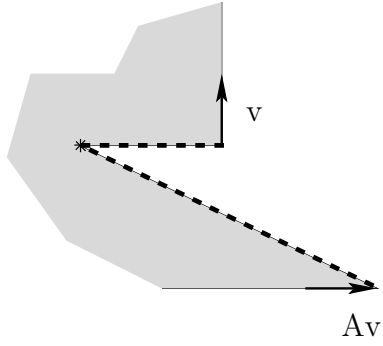}
\end{center}
\caption{Turning of tangent vectors to the boundary near a node.} 
\label{turning.fig}
\end{figure}

The essential fact is that $\theta_i\ge 0$ for each $i$.  
Indeed,
\be \label{rotation.eq}
\sin \theta_i=\frac{A_iv_i\times v_i}{\norm{A_iv_i}\norm{v_i}}
\ee
where $A_i$ is the affine monodromy around $b_i$.
(Here we are viewing the vector product on $\bR^3$ as
a scalar product on $\bR^2$.)
By direct calculation 
$Av\times v\ge0$ for any $A$ conjugate to $A_{(1,0)}$ 
and any vector $v\in\bR^2$.
Since $\theta_i$ measures the rotation of a vector under the linear map $A$,
$-\pi<\theta_i<\pi$, so $\sin \theta_i\ge 0$ implies $\theta_i\ge 0$.
Therefore $\int_D K_g = \sum_{i=1}^k \theta_i \ge 0$.
\end{proof}

\begin{rmk} To see that the weak inequality of Lemma~\ref{disk.lem} cannot
be replaced by a strict one, consider the
base diagram shown in Figure~\ref{blowup.fig}(b).  It represents a disk with 
one node such that any boundary compatible metric on the closure would have
total curvature equal to zero.
\end{rmk}

\begin{lemma} \label{Eulerchar.lem}
If \triple is the base of an \atf of a closed four-manifold
then $\chi(B)\ge0$.
\end{lemma}
Note that this lemma was observed by Zung~\cite{Zung.II} without proof.

\begin{proof}
In the arguments that follow we assume that the base is orientable;
if not, the double cover has an induced affine structure with twice as
many nodes.

If the base has no boundary then we can appeal to a theorem of 
Matsumoto~\cite{Matsumoto.ellipticsurf} on the structure of smooth Lefschetz
fibrations over surfaces of any genus, allowing us 
to find a flat metric on the
complement of a disk that contains all the nodes.
Specifically, Matsumoto showed that for any smooth Lefschetz fibration 
there is a presentation of the fundamental group of the base minus the
images of singular fibers such that each generator is
either a simple curve around the image of a singular fiber or
else is a curve along which the monodromy is trivial.
Recall that an almost toric fibration over a surface with no boundary is
smoothly equivalent to a Lefschetz fibration and the affine monodromy
is, up to taking a transpose, the same as the topological monodromy.
Therefore, the conclusion of Matsumoto's theorem holds for 
the affine monodromy of such an almost toric fibration.
Accordingly, choose a flat metric compatible with the affine structure
in a neighborhood of the generators of  $\pi_1(B-\cup_{i=1}^k b_i)$ along
which the monodromy is trivial.
Then, on the complement,
we have an affine disk with nodes that has, 
on a collar neighborhood of its boundary,
a flat metric compatible with the affine structure.
Adjusting the boundary of the disk if necessary so that it contains a 
linear segment, we can construct a boundary compatible metric on the disk.
Lemma~\ref{disk.lem} then implies that the base $B$ admits a metric
whose total curvature is non-negative, thereby forcing the Euler 
characteristic to be non-negative.

If the base has boundary we assume that either the genus of the 
base is at least one or the number of boundary components is at least two
(since the Euler characteristic of a disk is 
non-negative).
We work with a fundamental domain $(\wtilde D,\wtilde \cA)$ in the
universal cover $(\wtilde B,\wtilde \cA)$.

\begin{figure}
\begin{center}   	
	\psfragscanon
	\psfrag{g1}{$\gamma_1$}
	\psfrag{g2}{$\gamma_2$}
	\psfrag{g3}{$\gamma_3$}
	\includegraphics[width=2.5in]{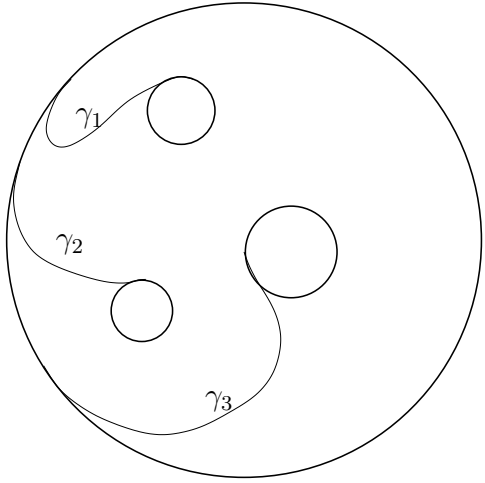}
\end{center}
\caption{Construction of fundamental domain for a thrice punctured disk.  }
\label{punctureddisk.fig}
\end{figure}

\begin{figure}
\begin{center}
	\psfrag{g1}{$\gamma_1$}
	\psfrag{g2}{$\gamma_2$}
	\psfrag{g3}{$\gamma_3$}
	\psfrag{g4}{$\gamma_4$}
	\psfrag{g5}{$\gamma_5$}
	\includegraphics[width=4in]{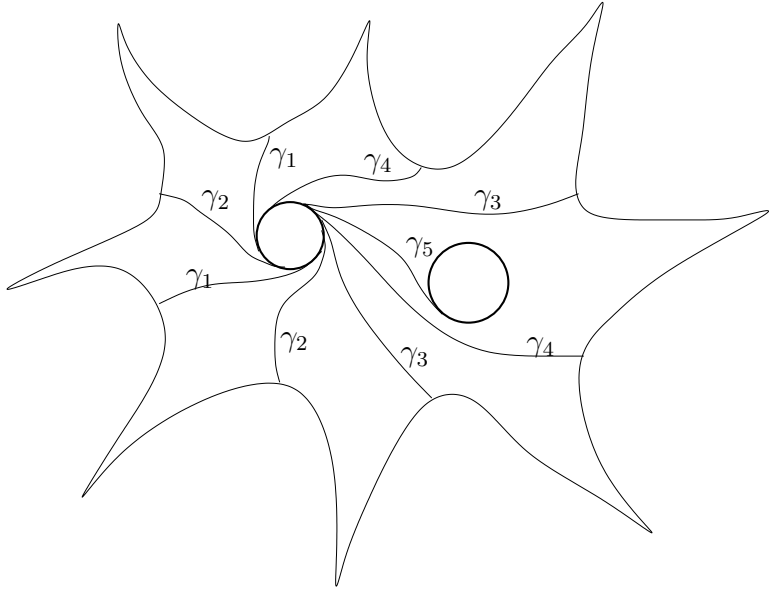}
\end{center}
\caption{ Construction of fundamental domain for a genus two surface with
two boundary components.}
\label{genusd.fig}
\end{figure}

Since the only angles that are well defined with respect
to an affine structure are multiples of $\pi$ we follow 
Benzecri~\cite{Benzecri.affine} and choose the fundamental 
domain so that the internal angle at any vertex
introduced in the universal cover is zero or $\pm \pi$.
For instance, to form the fundamental domain for a surface of genus zero with
four boundary components, we cut it open along three curves $\gamma_1$,
$\gamma_2$, and $\gamma_3$ as shown in Figure~\ref{punctureddisk.fig}.
Similarly, we form the fundamental domain for a surface of genus two with
two boundary components as suggested by Figure~\ref{genusd.fig}: 
for each $i$, identify the two smooth boundary components that intersect 
the curves $\gamma_i$, $i=1,\ldots, 4$ to get a genus two surface with two 
boundary components, then
cut it open along the curves $\gamma_i$, $i=1,\ldots, 5$.

The constraints on the topology of $B$ come from the geometry of
$(\wtilde D,\wtilde \cA)$.  Indeed, choosing a 
metric $\tilde g$ on $\wtilde D$ that is boundary compatible with
$\wtilde \cA$, the Gauss-Bonnet theorem implies
\be \label{constraint.ineq}
2\pi-\int_{\bdy\widetilde D} \kappa_{\tilde g}\, ds -\beta_{\tilde g} \ge 0
\ee
since $\int_{\wtilde D} K_{\tilde g} dA\ge 0$ by Lemma~\ref{disk.lem}. 

To estimate $\int_{\bdy\widetilde D} \kappa_{\tilde g}\, ds$, 
note that because $\bdy B$ (if
nonempty) is piecewise linear, all contributions
to the geodesic curvature along the smooth part of $\bdy\wtilde D$
come from pairs of arcs that each project to one arc in $B$.
Naming such a pair $\gamma_i, \gamma'_{i}$, their images in $\bR^2$ under
the developing map differ by an element of ${\rm Aff}(2,\bZ)$.
Specifically, for some orientation preserving element  
$\Psi\in{\rm Aff}(2,\bZ)$, $\gamma'_i=-\Psi(\gamma_i)$.
Since elements of ${\rm Aff}(2,\bZ)$ preserve the integer part of 
$\frac{\theta}{\pi}$ for any angle $\theta$, we have
\be 
\Bigl\lvert\int_{\gamma_i} 
\kappa_{\tilde g}\, ds+\int_{\gamma'_{i}} 
\kappa_{\tilde g}\, ds\Bigr\rvert<\pi
\ee 
for each such pair.

Let $d$ be the genus of the base $B$ and let $m$ be the number of
boundary components.
We now consider two cases separately: $\bdy B$ empty and $\bdy B$ non-empty.
Recall that we have assumed either $m\ge 1$ or $d\ge 2$. 

If $\bdy B$ is empty and then we can choose a fundamental domain in which
all but two of the vertices contribute $\pi$ to the total turning angle and
two of them contribute zero.  (This satisfies the requirement that the
sum of the internal angles of the fundamental domain must equal $2\pi$.)
Then 
\begin{align} 
\label{nobdyturning.eq}
\beta_{\tilde g} & = (4 d - 2)\pi\qquad {\rm and}\\
\label{nobdygeodcurv.eq}
\Bigl\lvert\int_{\bdy\widetilde D} 
\kappa_{\tilde g}\, ds\Bigr\rvert &< 2\pi d, 
\end{align}
in which case Inequality~\ref{constraint.ineq} implies $d=1$. 

Suppose now that $\bdy B$ is non-empty. 
The local convexity of $\bdy B$ implies that
the contributions to $\beta_{\tilde g}$ at any vertices of the fundamental
domain that project to vertices of $\bdy B$ are positive.
Meanwhile, the contributions
at the other vertices can be calculated exactly -- thanks to the turning
angle at each being a multiple of $\pi$.
Specifically, constructing the fundamental domain in analogy with 
Figures~\ref{punctureddisk.fig} and~\ref{genusd.fig}, we calculate
\begin{align} 
\label{bdyturning.eq}
\beta_{\tilde g} &\ge 2\pi(2d+m-1) \qquad {\rm and}\\
\label{bdygeodcurv.eq}
\Bigl\lvert\int_{\bdy\widetilde D} 
\kappa_{\tilde g}\, ds\Bigr\rvert &< \pi(2d+m-1).
\end{align}
These inequalities imply that either $d=0$ and $m=2$ or $d=1$ and $m=0$.

In all of the above cases $\chi(B)\ge 0$.
\end{proof}

To get more detailed information about the affine structure on the almost
toric bases in question we need 
the following standard fact about matrices that encode the monodromy
around a node: 
\begin{thm}(\cite{Moishezon.cxsurfaces}) \label{twelve.thm}
Suppose $\{A_i\}_{i=1}^k$ is a set of matrices in $SL(2,\bZ)$, each 
conjugate to 
$\left(\begin{smallmatrix}1&1\\0&1\end{smallmatrix}\right)$.
If $\Pi_{i=1}^k A_i=I$ 
then $k=12n$ and there is a finite sequence of
elementary transformations that yields the product $\Pi_{i=1}^k A'_i=I$ such
that $A'_i=
\left(\begin{smallmatrix}1&1\\0&1\end{smallmatrix}\right)$
if $i$ is even and 
$A'_i=
\left(\begin{smallmatrix}1 & 0\\-1&1\end{smallmatrix}\right)$
if $i$ is odd.
\end{thm}
Note that an elementary transformation on a cyclicly ordered set
of matrices is either
\be
\begin{split}
T_i:\{A_1,&\ldots, A_{i-1},A_i,A_{i+1},A_{i+2},\ldots ,A_k\}
\\ &\to
\{A_1,\ldots, A_{i-1},A_iA_{i+1}A_i^{-1},A_i,A_{i+2}, \ldots,A_k\}
\end{split}
\ee
for some $i$, or its inverse $T_i^{-1}$.
The relevance of the elementary transforms and Theorem~\ref{twelve.thm}
follows from:
\begin{obs} \label{elemtransf.obs}
Given an integral affine disk with nodes $(D,\cA)$, a choice of branch
curves $R=\cup_{i=1}^k R_i$ and an immersion 
$\Phi:(D-R,\cA)\ra(\bR^2,\cA_0)$
determines a representation of the affine monodromy in $SL(2,\bZ)$.
In particular, if the branch curves are indexed so that their intersections
with the boundary give an ordered set of points, say
$\{x_1,x_2,\ldots x_k\}$ agreeing with the orientation
of the boundary, then the monodromy along the boundary is $A_1A_2\ldots A_k$
where $A_i$ is the monodromy around the node $b_i$.
The elementary transformation $T_i$ then corresponds to replacing
the branch curve $R_i$ by a branch curve from $b_i$ to $x'_i$ where $x'_i$ 
is between $x_{i+1}$ and $x_{i+2}$.
Accordingly, 
we call a change in branch curves corresponding to $T_i$ or $T_i^{-1}$
an {\it elementary branch move}. 
\end{obs}

\begin{lemma} \label{diskcurvature.lem}
Let $(D,\cA)$ be an affine disk.  If
the monodromy around the boundary is trivial then there are $12n$ nodes
on the interior of $(D,\cA)$.
Furthermore, for any boundary compatible metric $g$,
\be \label{totalcurv.eq}
\int_D K_g \,dA = 2\pi n.
\ee 
\end{lemma} 

\begin{proof}

Theorem~\ref{twelve.thm} implies that the number of nodes
must be a multiple of $12$.
We proceed as in the proof of Lemma~\ref{disk.lem} choosing a set of
branch curves $R=\cup_{i=1}^k R_i$ and flat model $(\widehat D,\hat g)$.
However in this case, if there are nodes,
we appeal to Theorem~\ref{twelve.thm}, and 
choose the branch curves so that the monodromy
$A_i$ across each curve $R_i$ is
$\left(\begin{smallmatrix}1&1\\0&1\end{smallmatrix}\right)$
if $i$ is even and 
$\left(\begin{smallmatrix}1 & 0\\-1&1\end{smallmatrix}\right)$
if $i$ is odd.
A direct calculation then shows $\sum_{i=1}^{k} \theta_i=2\pi \frac{k}{12}
= 2\pi n$.
(It is easiest to check this on the vector $(1,0)$ but
it is independent of the choice of vector since the total monodromy
is trivial.)
As in the proof of Lemma~\ref{disk.lem} we have that 
$\int_D K_g=\sum_{i=1}^{k} \theta_i$ thereby proving the result.
\end{proof}

\begin{thm} \label{base.thm}
Suppose \triple is the base of an \atf of a closed four-manifold.
Then $(B,\cA)$ must be one of the following:
\bi
\item a disk with any number of nodes;
\item a cylinder or M\"oebius band with any number of nodes,
all of whose eigenlines are parallel to $\bdy B$ which is linear;
\item a closed surface with $12\chi(B)$ nodes.
\ei
\end{thm}

\begin{proof}
An easy way to construct an almost toric manifold whose base is a disk
with $k$ nodes and $v$ vertices is to start with the moment map image
of $\bC P^2\,\#\,(k+v-3)\overline{\bC P}^2$ (which is a polygon with 
$k+v$ vertices) and then perform $k$ nodal trades.

Suppose that $(B,\cA,\cS)$ is a sphere.  Cover $(B,\cA)$ with 
two disks, $(D_1,\cA)$ and $(D_2,\cA)$, such that 
$D_1\cap D_2=\bdy D_1=\bdy D_2$.  Assume 
that $D_2$ contains all of the nodes.
Then $(D_1,\cA)$ admits a flat metric compatible with $\cA$ and it extends
across $(D_2,\cA)$ as a boundary compatible metric $g$.
The Gauss-Bonnet theorem then implies that the total curvature on 
$D_2$ is $4\pi$, whereby Lemma~\ref{diskcurvature.lem} forces there to
be $24$ nodes.  (Consequently, if $B=\bR P^2$ then 
$(B,\cA,\cS)$ must have $12$ nodes.)  

If $B$ is a torus, a priori there could be nontrivial 
monodromy along generators
of $\pi_1(B)$ that is balanced by monodromy around nodes $\cup_{i=1}^k b_i$.
This possibility however is precluded by the theorem of 
Matsumoto referred to in the proof of Lemma~\ref{Eulerchar.lem}.
Therefore we can proceed as for a sphere, constructing 
a compatible flat metric on the
complement of the interior of a disk that contains all the nodes.
Then $\chi(B)=0$ implies that the curvature on that disk is zero, and
hence by Lemma~\ref{diskcurvature.lem} the disk must contain no nodes.

Now assume that $B$ is a cylinder and that there are no vertices on the
boundary, performing nodal trades if necessary.
Example~\ref{sphereovertorus.ex} 
shows that there is no bound on the number of 
nodes.  However, the monodromy around the nodes is quite restricted: 
Lemma~\ref{cylnodes.lem} below implies
that we can choose a monodromy representation in which the monodromy around
all the nodes the same, i.e. all the nodes belong to a disk in which
there is a well defined line -- the eigenline through the nodes.
It remains to check that this line is parallel to the components of 
$\bdy B$ which in turn are parallel to each other.
(Recall that the boundary is geodesic 
because the total space of the fibration is a closed manifold.)

For the monodromy presentation choose a base point $b\in B-\cup_{i=1}^k b_i$, 
a basis for $T_b B$, 
and simple loops $\gamma_i$, $1\le i\le k$, based at $b$ such that
\ben
\item each $\gamma_i$ winds around exactly one node $b_i$, positively oriented
as the boundary of the disk it bounds and
\item the affine monodromy along each $\gamma_i$ is $A_{(1,0)}$.
\een
Now let $\eta_1$, $\eta_2$ be loops based at $b$ that generate 
$\pi_1(B)$ and are such that 
\be \label{monodromy.eq}
\eta_2^{-1}\eta_1=\gamma_k\cdots \gamma_1
\ee
so that there are
no nodes between $\eta_1$ and one component of $\bdy B$ and no nodes
between $\eta_2$ and the other component of $\bdy B$.
Since the affine monodromy along each $\gamma_i$ is $A_{(1,0)}$, 
Equation~\ref{monodromy.eq} implies that the affine
monodromy along $\eta_2^{-1}\eta_1$ is $A_{(1,0)}^k.$  
It is easy to check that this can happen only if the monodromy along
both $\eta_1$ and $\eta_2$ is a power of $A_{(1,0)}$ whenever $k\ne 0$.  
Geometrically, this means that the boundary components are parallel to
the common eigenline of the nodes.  Note that it also implies $\bdy B$
cannot have any vertices.

If there are no nodes (the case $k=0$) then the monodromy along $\eta_1$ is 
the same as along
$\eta_2$.  Furthermore, this monodromy must be conjugate to a power of 
$A_{(1,0)}$ since a tangent vector to either linear boundary component 
must be invariant under the monodromy.  (This is most easily seen in
a base diagram of a fundamental domain.)  Therefore the boundary
components are parallel to each other.

\end{proof}

\begin{lemma} \label{cylnodes.lem}
Suppose $(B,\cA)$ is a cylinder whose boundary is linear with respect to $
\cA$.
Assume that there are nodes in the affine structure $\cA$.
Then there is a presentation of the affine monodromy such that the monodromy
around each node is $A_{(1,0)}$, i.e. all the eigenlines are parallel.
\end{lemma}

\begin{proof}
Choose a fundamental domain 
$(\wtilde D,\tilde \cA)\subset(\wtilde B,\tilde \cA)$ 
and a boundary compatible metric $\tilde g$.
Then Equations~\ref{bdyturning.eq} and~\ref{bdygeodcurv.eq} imply, with
$d=0$ and $m=2$, that 
$\beta_{\tilde g}\ge 2\pi$ and  
$\bigl\lvert\int_{\gamma} \kappa_{\tilde g}\, ds \bigr\rvert<\pi$.
With these constraints the Gauss-Bonnet theorem forces 
$\int_{\wtilde D} K_{\wtilde g}\ dA<\pi$.
Our proof by contradiction amounts to construction of a boundary
compatible metric that violates this bound.    

If there is only one node, then the conclusion of the lemma is trivial.
Therefore, let $\{b_i\}_{i=1}^k$, $k\ge 2$, be the nodes indexed so  
that the monodromy around $b_2$ is not the same as around $b_1$.
Assume without loss of generality that $\bdy\wtilde D$ 
contains at least four linear segments, one pair
of which (say $L_1,L_2$) determines the other.
Choose a set of branch curves $R=\cup_{i=1}^k R_i$ such that 
$R_i\cap \bdy \wtilde D\in L_i$ for
$i=1,2$ and an immersion $\Phi:(\wtilde D-R,\cA)\ra(\bR^2,\cA_0)$ such
that the monodromy around $b_1$ is $A_{(1,0)}$.
Then following the proof of Lemma~\ref{bdycompat.lem},
we construct a boundary compatible metric $\tilde g$ on 
$(\widetilde D,\widetilde A)$ whose total curvature is 
$\int K_{\tilde g}\, dA= \sum_{i=1}^k \theta_i$ where 
$\theta_i$ measures the curvature contribution of $b_i$ as in
the proof of Lemma~\ref{disk.lem}.
Since $\theta_i\in [0,\pi)$ for each $i$ (as shown in that proof), we only
need to verify that we can choose the fundamental domain 
so that $\theta_1+\theta_2\ge\pi$. 

Referring to Figure~\ref{turning.fig} and suppressing indices, we note that 
\be
\tan \theta = \frac{Av\times v}{Av\cdot v}
\ee
where $A$ is the monodromy around a node.
With $A=A_{(p,q)}$ and $v=(x,y)$ we calculate
\begin{align}
Av\times v&=(qx-py)^2 \label{cross.eq}\\
Av\cdot v&= (1-pq)x^2+(p^2-q^2)xy+(1+pq)y^2. \label{dot.eq}
\end{align}

Allowing the affine lengths of 
$L_1,L_2$ to be sufficiently small, we can choose the 
fundamental domain so that the vectors $v_1,v_2$ have any direction we want.
Recalling that we have chosen $A_1=A_{(1,0)}$, 
Equations~\ref{cross.eq} and~\ref{dot.eq} imply $\theta_1$ is
maximized by taking $v_1=(-1,2)$ in which case $\tan \theta_1=4/3$.  
For $A_2=A_{(p,q)}$ we can assume $q\ge 1$ 
since $q=0$ would make the eigenvectors parallel 
and the vector $(p,q)$ is defined only up to sign.
Furthermore, we can choose (without loss of generality) $p$ to be 
positive and arbitrarily large since it is defined only mod $q$.
(To change the value of $p$ we can modify our choice of affine immersion
$\Phi$ in a way that causes $A_2$ to be conjugated by a
power of $A_{(1,0)}$ -- and therefore leaves $A_1$ unchanged.)

Having chosen $v_1=(-1,2)$ the bound $\theta_1+\theta_2<\pi$ will violated if
\begin{align}
\tan\theta_2 &<0 \ \ {\rm and}\\
\abs{\tan\theta_2}\le \tan\theta_1 & =\frac{4}{3}.
\end{align}
We do this by choosing $v_2=(p+2,q)$ so that
\begin{align}
A_2v_2\times v_2 & =4q^2 \ \ {\rm and}\\
A_2v_2\cdot v_2 & = (1-2q) p^2 + 4(1-q) p + (4+q^2-2q^3).
\end{align}
Choosing $p$ large enough we have that 
$\tan \theta_2$ is negative and as 
close to zero as we like.
\end{proof}

\subsection{Disk base}

The goal of this section is to prove that if $\pi:(M,\omega)\ra(B,\cA)$ is
a closed almost toric four-manifold fibering over a disk, then $M$ is 
diffeomorphic to a toric manifold.  In other words,

\begin{thm}  \label{attotoric.thm}
Suppose $\pi:(M,\omega)\to(D,\cA,\cS)$ 
is a closed almost toric manifold whose base is a disk.  
Then there is a symplectic structure $\omega'$, 
deformation equivalent to $\omega$, such that
$(M,\omega')$ admits a toric fibration.
\end{thm}

The classification of toric manifolds immediately implies
\begin{cor} \label{attotoric.cor}
If $\pi:(M,\omega)\to(D,\cA,\cS)$ 
is a closed almost toric manifold whose base is a disk, then
$M$ is determined up to diffeomorphism by the integer $k+n$ where
$k$ is the number of nodes and $n$ is the number of vertices, unless
$k+n=4$ in which case $M$ can have one of two topological types.
\end{cor}

The essential idea of the proof of Theorem~\ref{attotoric.thm}
is simple: just slide the nodes on the
interior of $B$ along eigenrays until they hit the boundary,
thereby performing the inverse of a nodal trade.
However, the almost toric base in question may be such that:
\ben
\item there is no set of disjoint eigenrays connecting the nodes to
the boundary (along which to slide the nodes) and/or 
\item sliding a node all the way to the boundary might 
produce a change in topology by creating an
orbifold singular point.
\een

The first issue is addressed by Lemma~\ref{vectors.lem} which allows us
to assume, since we are only interested in the topology of the total
space, that all nodes are sufficiently close to the boundary that
there is a \lq\lq good\rq\rq\ set of branch curves.
Specifically, the branch curves can be chosen so that each one belongs either
to an eigenline or to a neighborhood of a boundary point in which 
the eigenline through the node is parallel to the base.  
The influence of nodes of the latter type is explained by the
discussion of Section~\ref{atblowup.sec} which 
shows that the node is the result of an almost toric blow-up.

We resolve the second issue by providing 
(in the proof of Lemma~\ref{goodorder.lem})
an algorithm to appropriately modify an almost toric disk base
without changing the topology of the total space it defines.

For simplicity of exposition, we assume that there are no vertices
on the boundary  (so the zero-stratum of the base is empty).
Otherwise we start by trading all vertices for nodes.

To begin, we indicate a set of the data on a base diagram that defines the 
topology of the total space.
\begin{lemma} \label{vectors.lem}
Suppose $\pi:(M,\omega)\to(D,\cA,\cS)$ 
is a closed almost toric manifold whose base is a disk with no vertices
on the boundary.
Suppose $V=\Phi(D-R)\subset (\bR^2,\cA_0)$ is the domain of a base 
diagram,
where $R=\cup_{i=1}^k R_i$ is a set of branch curves.
Let $u_i$, $i=1,\ldots k$ be the inward-pointing primitive integral vectors 
normal to the connected components of $\bdy V$,
indexed so that they rotate non-negatively (counterclockwise).
Then the set of vectors $\{u_1,\ldots u_k\}$, up to
cyclic permutation and the action of $GL(n,\bZ)$, 
determines the diffeomorphism type of $M$.  
\end{lemma}
Note that $\Phi$ is an embedding since $\bdy D$ is locally convex.
Also, each vector $u_i$ really should be viewed as a
covector defining the corresponding connected component of $\bdy V$.

\begin{remark} \label{fan.rmk}
For the reader familiar with complex algebraic toric varieties, 
the vectors $u_i$ define a complete fan.
The toric variety defined by the fan will in general have
orbifold singularities and hence not be diffeomorphic to $M$.
It will fail to be smooth precisely when $u_i\times u_{i+1}>1$
for some $i$.   
\end{remark}

\begin{proof}
Because $(D,\cA)$ has no vertices its boundary is geodesic with
respect to $\cA$.
Therefore $M$ is the boundary reduction (along one line)
of a symplectic manifold  $(M',\omega')$ that is 
a smooth Lefschetz fibration over a disk.  
As such, the diffeomorphism type of $M'$ is
completely determined by the monodromy. 
The only additional data needed to determine $M$ is the homology class
of a regular fiber that gets collapsed during the boundary reduction
(the {\it collapsing class}).
This homology class is well defined only with respect to an
arc in the base that runs from the image of the
regular fiber to the boundary.

The ordered set of vectors $\{u_1,\ldots,u_k\}$  defines 
the monodromy because
for each pair $\{u_i,u_{i+1}\}$ there
is a unique matrix $A_i$ conjugate to $A_{(1,0)}$ such that $A_iu_i=u_{i-1}$
(mod $k$).
Furthermore, the collapsing class can be defined by $u_i$ for any $i$ as
follows: Choose an embedded arc $\gamma$ connecting a point $\Phi(b)$ on the
interior of $V$ and a point in the connected component of $\bdy V$ defined
by $u_i$.  Viewing $u_i$ as a covector in $T^*\bR^2$ and pulling
back via $\Phi$ we get a covector in $T^*_b D$.  
The collapsing class with respect to $\Phi^{-1}(\gamma)$ is then the element
of $H_1(F_b,\bZ)$ represented by integral curves of the vector
field $X$ such that $\omega(X,\cdot)=\Phi^*u_i$.  (Here $F_b$ is the regular
torus fiber over the point $b\in D$ and $\omega$ is the symplectic structure
defined by the almost toric fibration.)

Since cyclicly permuting the vectors  $u_1,\ldots,u_k$ 
has no effect on the monodromy presentation or the collapsing class
it also has no effect on the topology.
Furthermore, changing the vectors by applying an element of  
$GL(2,\bZ)$ amounts only to changing the isomorphism between 
$H_1(F_b,\bZ)$ and $\bZ^2$.
\end{proof}

\begin{rmk}
Unless the monodromy along the boundary of the base (with no vertices)
is trivial, the collapsing class is determined by the monodromy.  Indeed,
$A_1A_2\cdots A_k u_k=u_k$ and, unless $A_1A_2\cdots A_k=I$, the primitive
integral vector $u_k$ is unique up to sign.
In contrast, if $A_1A_2\cdots A_k=I$ then any vector could determine a
collapsing class with respect to a fixed arc.  
Upon proving Theorem~\ref{attotoric.thm}, Corollary~\ref{attotoric.cor}
implies the diffeomorphism type of the total space is
independent of this choice of vector.
Since the boundary is geodesic, the proof of Lemma~\ref{diskcurvature.lem} 
forces the number of nodes to be $12$ and therefore 
the total space is diffeomorphic to the elliptic surface $E(1)$.
Accordingly, the possibility of choosing any vector to determine the 
collapsing class whenever $A_1A_2\cdots A_k=I$  
reflects the very large diffeomorphism group of $E(1)$.
\end{rmk} 

A natural question is what sequences of vectors $\{u_1,\ldots, u_k\}$
can be the normal vectors to the connected components of $\bdy V$ 
where $V=\Phi(D-R)$ as above?
The primary constraint is that 
\be
A_iu_i=u_{i-1}
\ee
for each $i$ (mod $k$)
where $A_i$ is some matrix conjugate to $A_{(1,0)}$.
The action of $A_i$ on any vector $v$ can be rewritten in terms of its
eigenvector $e_i$ as:
\be
A_i v= v - (v\times e_i) e_i \label{Aiv.eq}
\ee
Therefore, the constraint can be rewritten as 
\be
u_{i-1}=u_i-(u_i\times e_i) e_i
\ee
for some primitive integral vector $e_i$.  
The only other constraint on $\{u_1,\ldots u_k\}$ is
that the vectors rotate exactly once around the origin.
This motivates

\begin{defn} \label{defining.def}
An ordered set of primitive integral vectors $(u_1,\ldots, u_k)$ 
is a {\it defining set} for a closed almost toric manifold fibering
over a disk if for each $i$ there is an integer $n_i$ and primitive integral
vector $e_i$ such that 
\ben
\item $u_i\times e_i=n_i$,
\item $u_i-u_{i-1}=n_i e_i$ (mod $k$),
\item  $u_j\ne u_1$ for some $j\ne 1$ and
\item if $u_m= u_1$ for some $m$, then either $u_i=u_1$ for all $i\le m$ or
else $u_i=u_1$ for all $i\ge m$.
\een
\end{defn}
Note that for any defining set $\{u_1,\ldots u_k\}$ the corresponding
set of integers $\{n_1,\ldots n_k\}$ is such that 
$n_i^2=u_{i-1}\times u_i$ for each $i$ (mod $k$).
Furthermore, the definition of $e_i$ is such that $n_i\ge 0$ for all $i$.

Following up on Remark~\ref{fan.rmk}, the fan defined by $\{u_1,\ldots u_k\}$ 
is the fan for a smooth variety if and only if $n_i=1$ for all $i$.
Accordingly, being able to slide a node $b_i$ into the boundary of an almost
toric base without a change in topology of the total
space requires $n_i=1$.
If $n_i=0$ the eigenline through $b_i$ is the eigenvector of $(A_i^{-1})^T$
and hence is parallel to the boundary
inside a disk containing $R_i$, the branch curve with one endpoint at
$b_i$.  As mentioned before, in this case the node is the result of
blowing up.

\begin{lemma} \label{goodorder.lem}
Suppose $\{u_1,\ldots, u_k\}$ is a defining set for a closed
symplectic four-manifold $(M,\omega)$.
There is a sequence of elementary branch
moves (defined in Observation~\ref{elemtransf.obs})
that yields a new defining set $\{u'_1,\ldots u'_k\}$ such that
$u'_{i-1}\times u'_i=(n'_i)^2\in\{0,1\}$ for all $i$.
\end{lemma}

\begin{proof}

The defining set determines a corresponding set of monodromy matrices
$\{A_1,A_2,\ldots A_k\}$.
The elementary branch move corresponding to the elementary transformation
$T_j$ causes $u_j$ to be replaced by $A_ju_{j+1}$, leaving
the other vectors $u_i$, $i\ne j$, unchanged. 

Let $\tau_j$ be the induced action on the integers $\{n_1,n_2,\ldots n_k\}$.
Then 
\be \tau_j(n_i)=n_i \ \ {\rm if \ } i\ne j, j+1, \ee
\begin{align}
\tau_j(n_j)& = \sqrt{u_{j-1}\times A_j u_{j+1}} \\
           & = \sqrt{A_ju_j\times A_j u_{j+1}} \\
           & = \sqrt{u_j\times u_{j+1}} = n_{j+1} 
\end{align}
and
\begin{align}
\tau_j(n_{j+1}) & = \sqrt{A_ju_{j+1}\times u_{j+1}} \\
           & = \sqrt{\left(u_{j+1}-\left(u_{j+1}\times e_j\right) e_j\right)
                 \times u_{j+1}} \\
           & = \sqrt{(u_{j+1}\times e_j)^2} = \lvert u_{j+1}\times e_j\rvert 
\end{align}
where the second equality follows from Equation~\ref{Aiv.eq}.
Therefore, performing a sequence of 
elementary branch moves corresponding to $T_{j+m}\cdots T_{j+1}T_j$
has the effect,
via $\tau_{j+m}\cdots\tau_{j+1}\tau_j$, of removing $n_j$ from the set and
inserting $\lvert u_{j+m+1}\times e_j\rvert$.

Assume without loss of generality that 
$n_1=N\ge n_i$ for all $i$.
Also assume that $N\ge 2$ (for otherwise
our initial sequence would satisfy the conclusion of the lemma).

If $\abs{u_j\times e_1}<N$ for some $j$, then performing a branch move
that corresponds to $T_{j-1}\cdots T_2T_1$
removes $n_1=N$ and replaces it by a strictly smaller non-negative integer.
Redefining $u_1$ we could apply the same argument repeatedly.
Therefore the
only obstruction to achieving $n_i\in\{0,1\}$ for all
$i$ would be if at some stage 
$\abs{u_i\times e_1}\ge N$ for all $i$.  Assume this is true.

\begin{figure}
\begin{center}
        \psfrag{2N}{$2N$}
        \psfrag{f1}{$f_1$}
        \psfrag{uk}{$u_k$}
        \psfrag{ul}{$u_l$}
        \psfrag{u1}{$u_1$}
        \psfrag{ujm1}{$u_{j-1}$}
        \psfrag{ulm1}{$u_{l-1}$}
        \psfrag{uj}{$u_{j}$}
        \psfrag{n1e1}{$n_1e_1$}
        \psfrag{njej}{$n_je_j$}
        \psfrag{nlel}{$n_le_l$}
	\includegraphics[width=3in]{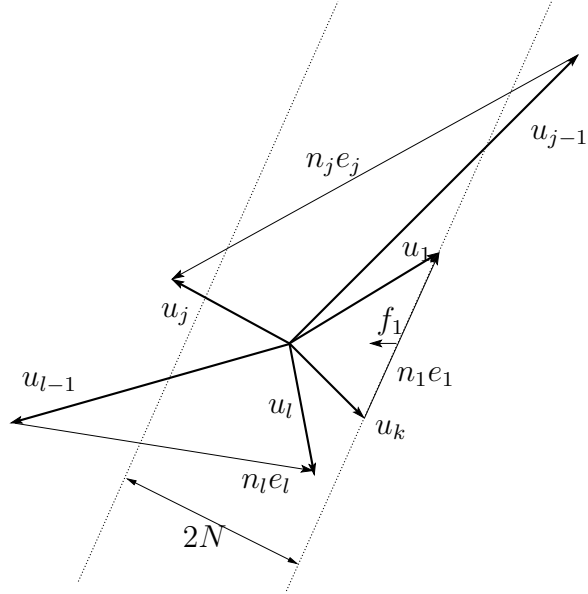}
\end{center}
\caption{Vectors of a defining set.}
\label{vectors.fig}
\end{figure}

Let $f_1$ be a primitive integral vector such that $e_1\times f_1=1$.
Then we can write any vector $v$ as a linear combination of $e_1,f_1$
where $e_1\times v$ gives the $f_1$ coefficient and 
$v\times f_1$ gives the $e_1$ coefficient.

Since the equation $e_1\times f_1=1$ defines $f_1$ modulo
an integer multiple of $e_1$, we can choose
$f_1$ so that $u_k\times f_1 < 0$ and 
and $u_1 \times f_1 \ge 0$.  (See Figure~\ref{vectors.fig}.)  
Furthermore, we can strengthen the last inequality to 
$u_1\times f_1 > 0$ because $u_1\times f_1=0$ and the primitivity of 
$u_1$ would imply
$u_1=-f_1$ whence $N=\sqrt{u_1\times e_1}=\sqrt{e_1\times f_1}=1$
contradicting the assumption $N\ge 2$.

Because the angle between $u_{i-1}$ and $u_i$ is less that $\pi$ for any $i$, 
there must be some minimal $j$ such that 
$u_j\times e_1\le -N$.  
Writing
$e_{j}= x e_1 + y f_1$ where $x,y$ are both integers, our choice
of $f_1$ implies $y>0$.
Meanwhile,
\begin{align}
n_j & = u_j\times e_j=u_{j-1}\times e_j \\
    & = x (u_{j-1}\times e_1) + y (u_{j-1}\times f_1),
\end{align}
but $u_{j-1}\times e_1\ge N$ and $u_{j-1}\times f_1> 0$ so the only way
to have $n_j\le N$ is if $x<0$.
Furthermore,
\begin{align}
n_j & = u_j\times e_j \\
    & = x (u_j\times e_1) + y (u_j\times f_1)\\
\end{align}
where $u_j\times e_1\le -N$.  With $y>0$ and $x<0$, the constraint that
$n_j\le N$ forces $u_j\times f_1\le 0$.
Since $u_j\times f_1=0$ would imply $u_j=e_1$ and thereby $N=1$, we
find $u_j\times f_1< 0$, i.e. the $e_1$ component of $u_j$ must be
negative.

By symmetry, the same argument for $u_{l-1},u_{l}$, 
where $l$ is the maximal index for which $e_1\times u_{l-1}\ge N$, would show 
that the $e_1$ component of $u_{l-1}$ must be positive.

Since both of these conditions of the $e_1$ components of $u_j$ and $u_{l-1}$
cannot be met, the assumption that $\lvert u_i\times e_1\rvert\ge N\ge 2$
must have been false.
\end{proof}

\begin{proof}[Proof of Theorem~\ref{attotoric.thm}]
Assume without loss of generality that the almost toric fibration is
over a disk that has no vertices.
Let $\{u_1,\ldots u_k\}$ be a defining set of vectors 
arising from a particular base diagram.
Invoking Lemma~\ref{goodorder.lem} we can, by varying the base diagram
without changing the fibration,
find a new defining set $\{u'_1,\ldots, u'_k\}$ such that
$n_i^2=u'_{i-1}\times u'_{i}\in\{0,1\}$ for all $i$.  In this base diagram
the branch curves need not be linear.  

Now, allowing the fibration and symplectic structure to vary,  
we construct a new almost toric base
$(D,\cA',\cS')$ that defines the same smooth four-manifold but has
a more amenable base diagram.
Indeed, letting $l\le k$ be the number of distinct vectors
in the defining set, we construct the base diagram as follows:
\ben
\item Choose a convex polygon such that $\{u'_1,\ldots, u'_k\}$ is a set
of inward pointing normal vectors that rotate non-negatively counterclockwise.
This will be a polygon with $l$ sides.
\item For each $i$ such that $u'_{i}\ne u'_{i-1}$ (mod $k$),
place a dotted line
segment $\eta'_i$ in the polygon so that it has one endpoint at 
the vertex between the sides with normal vectors $u'_{i-1},u'_i$ 
and has $e'_i$ as a tangent vector.  Do this so that the $\eta'_i$ are
all disjoint.  
After placing an asterisk at the interior endpoint of each $\eta'_i$,
this will be the base diagram for an almost toric manifold.
\item For each $j$ such that $u'_j=u'_{j-1}$, perform an
almost toric blowup on the edge defined by $u_j$.
\een

By construction we can now slide all nodes to the boundary, perform
$k-l$ almost toric blowdowns, and then perform $k-l$ toric blowups.
The result will be a toric fibration of $(M,\omega')$ for some
symplectic structure $\omega'$.

To see that $\omega'$ and $\omega$ must be deformation equivalent note that
one could interpolate between the initial base diagram and the final base
diagram via a one parameter family of base diagrams.
Accordingly one can find a one parameter family of fibrations interpolating
between the initial and final ones.
\end{proof}

\subsection{Diffeomorphism classification} \label{class.sec}

Our main theorem now follows easily from Theorems~\ref{base.thm}
and~\ref{attotoric.thm} and the catalog of examples given in 
Section~\ref{closed.sec}.

\begin{proof}[Proof of Theorem~\ref{diffclass.thm}]
We know the base must have non-negative Euler characteristic.

If the base is a disk then Theorem~\ref{attotoric.thm} implies that it must
be a rational surface.

If the base is a cylinder or Mo\"ebius band, Theorem~\ref{base.thm} implies
that any nodes in the base must be the result of an almost-toric blow-up
(since the eigenlines must be parallel to the boundary).
Blowing down, we get a base with no nodes whose monodromy is 
$\left(\begin{smallmatrix}1&n\\0&1\end{smallmatrix}\right)$.
The parity of $n$ determines whether the total space is a trivial or
non-trivial sphere bundle over $T^2$.
Of course, we are free to blow up an arbitrary number of times so long
as the blow-ups are small enough.

If the base is a sphere, $\bR P^2$, torus or Klein bottle,
then Theorem~\ref{base.thm} implies there must be $24$ or $12$ nodal fibers in 
the first two cases respectively, and $0$ nodal fibers in either of
the last two cases.
As a smooth fibration, an almost toric fibration over $S^2$ is equivalent to
a genus-one Lefschetz fibration with $24$ singular fibers, which is
well-known to be diffeomorphic to a K3 surface.
When there are $12$ nodal fibers we proved that the total space is
diffeomorphic to $E(1)_{2,2}$ which is known to be the Enriques surface.
When the base is a torus, since there are no singular fibers the list
of examples from Section~\ref{closed.sec} is complete.
Similarly, the case of a Klein bottle base is covered by 
Example~\ref{torioverKlein.ex} since there are no nodes in the base. 

\end{proof}

\subsection{Other classifications}

{\sc Fiber-preserving symplectomorphism:}
Two Lagrangian fibrations are equivalent if and only if there is
a fiber-preserving symplectomorphism between them.  The problem of
determining the data required to specify a Lagrangian fibration up to
fiber-preserving symplectomorphism has been studied in several cases:
\ben
\item Duistermaat~\cite{Duistermaat.actionangle} 
solved this problem for regular Lagrangian fibrations: one
needs the affine structure of the base and a Lagrangian Chern class that 
measures the obstruction to the existence of a Lagrangian section.
\item Boucetta and Molino~\cite{BoucettaMolino.fibrations} solved this
problem for locally toric fibrations.  
The data consists, in our language, of the base
$(B,\cA,\cS)$ and a generalization of Duistermaat's
Lagrangian Chern class.  They also determined that any choice
of Chern class and base $(B,\cA,\cS)$ (with the correct local structure) 
can be realized by such a fibration.
\item Zung~\cite{Zung.II} made a significant generalization to  fibrations 
with a class of singularities he calls 
\lq\lq non-degenerate topologically stable\rq\rq.
The data includes the base $(B,\cA,\cS)$, the fiber-preserving 
symplectomorphism type of the neighborhood of each singular fiber, 
some global topological data and an appropriately generalized Lagrangian
Chern class.
The work of Vu Ngoc, S.~\cite{VuNgoc.focusfocus} 
shows that the structure of the fibration
near a singular fiber is delicate information already in dimension four where
he found a Fourier series type invariant for the neighborhood 
of a focus-focus (nodal) singularity.
\een

For closed manifolds one could hope for a complete classification that
specifies what fibrations can occur.
\ben
\item For regular Lagrangian fibrations of four-manifolds, this problem
was completely solved by Mishachev~\cite{Mishachev.Lagbundles}.
\item For toric fibrations this amounts to Delzant's 
theorem~\cite{Delzant.moment} and the classification of polytopes
satisfying the appropriate integrality conditions at each vertex.
\item To extend to locally toric fibrations in dimension four one needs to
treat the cases when the base is a cylinder, M\"oebius band or Klein bottle.
The first two cases amount to an easy exercise 
since the Lagrangian Chern class vanishes and the possible bases 
$(B,\cA,\cS)$ are easy to specify.
Meanwhile, Mishachev's work~\cite{Mishachev.Lagbundles} significantly informs
the case when the base is a Klein bottle.
\item Extending to the almost toric case would require an understanding of
all the affine structures with nodes that can occur on $S^2$.  This question
is of independent interest in the context of mirror symmetry 
(cf.~\cite{GrossSiebert.logdegen}).
\een

{\sc Global symplectomorphism:}  If one has a classification up to
fiber-preserving symplectomorphism this amounts to deciding which
fibrations are equivalent via a global symplectomorphism.  Even in
the case of closed toric manifolds this is a nontrivial problem.
Mishachev conjectured that two Lagrangian torus bundles over a torus
are symplectomorphic if and only if they fiber-preserving symplectomorphic.
This could also be an interesting question for almost toric fibrations of
the K3 (and hence Enriques) surface.

{\sc Global diffeomorphism:}  This paper solves this problem for closed
almost toric four manifolds.  Smith~\cite{Smith.torusfibr} considered
and solved this problem for Lagrangian fibrations that are locally Lefschetz
and are such that a regular fiber is non-trivial in homology (thereby
excluding the torus bundles over tori that have $b_1=2$).
One could hope to carry out the program of this paper in higher dimensions,
but already in dimension six the possibilities for bases becomes quite vast.

{\sc Weak deformation:}  If two Lagrangian
fibrations are known to have diffeomorphic total spaces, one can ask whether
the pull back (via some diffeomorphism) of one symplectic structure can be 
connected to the other symplectic structure via a path of symplectic 
structures, i.e. whether they are weakly deformation equivalent.
One way to verify such a relationship is to connect the two
fibrations by a path of fibrations.  We conjecture, for
instance, that the symplectic structures on any pair of almost toric 
K3 surfaces are weakly deformation equivalent.

\end{document}